\documentclass{amsart}

\usepackage{epsfig}
\usepackage{amssymb}
\usepackage{graphicx}
\usepackage{epstopdf}
\usepackage{amscd}
\usepackage{graphics}
\newtheorem{Theorem}{Theorem}[section]
\newtheorem{Proposition}{Proposition}[section]
\newtheorem{Corollary}{Corollary}[section]

\def\proof{\par{\it Proof}. \ignorespaces}
\def\endproof{{\ \vbox{\hrule\hbox{%
     \vrule height1.3ex\hskip0.8ex\vrule}\hrule }}\par}
\newenvironment{Proof}{\proof}{\endproof}

\theoremstyle{definition}
\newtheorem{Definition}[Theorem]{Definition}
\newtheorem{Example}[Theorem]{Example}
\newtheorem{Notation}[Theorem]{Notation}
\newtheorem{Conjecture}[Theorem]{Conjecture}

\theoremstyle{remark}
\newtheorem{Remark}[Theorem]{Remark}

\numberwithin{equation}{section}

\begin{document}

\title{Toda lattice, cohomology of compact Lie groups and finite Chevalley groups}

\author{Luis Casian}
\address{Department of Mathematics, Ohio State University, Columbus,
OH 43210}
\email{casian@math.ohio-state.edu}
\author{Yuji Kodama$^*$}
\thanks{$^*$Partially
supported by NSF grant DMS0404931}

\address{Department of Mathematics, Ohio State University, Columbus,
OH 43210}
\email{kodama@math.ohio-state.edu}

\keywords{Toda lattice, cohomology of compact group, finite Chevalley groups, flag manifold,
Painlev\'e divisor, Frobenius map, Lefschetz fixed point theorem}

\maketitle


\thispagestyle{empty}
\pagenumbering{roman}
\tableofcontents

\pagenumbering{arabic}
\setcounter{page}{1}

\section{Introduction}
In this article we establish  a connection between the  Toda lattice defined for a real
split semisimple Lie algebra $\mathfrak g$ and the integral 
cohomology of a real flag manifold.   We denote by  $\check {\mathfrak g}$ the  real split Lie algebra with Cartan matrix given by  the transpose of the Cartan matrix of $\mathfrak g$.  Note that ${\mathfrak g}=\check {\mathfrak g}$ if ${\mathfrak g}$ is simple and not of type $B$ or $C$.  

The main results establish a \lq\lq partial dictionary\rq\rq   between objects related to the Toda lattice and objects associated to the flag manifold.  Then we have on one side of this dictionary 
an integrable system, the Toda lattice associated to $ {\mathfrak g}$ and related  objects,  in particular, the set of singularities
(blow-up points) of the trajectories of the Toda lattice and the algebraic varieties of
the set of zeros of Schur polynomials giving the solutions of the
nilpotent Toda lattice. On the other side of the dictionary we have
 the {\it real} flag manifold $\check G/\check B$  of the Lie group $\check G$ with $\check{\mathfrak g}={\rm Lie}(\check G)$, $\check B$ a Borel subgroup of $\check G$,
 and other Lie-theoretic objects related 
to it, such as a maximal compact Lie subgroup $\check K$ of $\check G$ and the rational cohomology
of $\check K$. In this setting of real split semisimple Lie groups,  the flag manifold is $\check G/\check B=\check K/\check T$ where $\check T$ is a {\it finite} group.
Thus  $H^*(\check G/\check B; {\mathbb Q}) = H^*(\check K; {\mathbb Q})$  although $\check G/\check B$ and $\check K$ have very different integral cohomology.  Therefore one additional  object which is connected to the real flag manifold
is the cohomology ring $H^*(\check K; {\mathbb Q})$.  In this paper we then derive  certain surprising  connections among  the Toda lattice, a maximal compact Lie subgroup
$\check K$ of $\check G$  and a finite Chevalley group $\check K({\mathbb F}_q)$
over a field ${\mathbb F}_q$ with $q$ elements. These relations amount to  a description of the rational cohomology  $H^*(\check K;{\mathbb Q})$ in terms of the blow-up points along trajectories of the Toda flow. Moreover this can be refined to  give a description of  the  integral cohomology of $\check G/\check B$ in terms of the Toda lattice.  
\smallskip

It is  possible to see the relation between the Toda lattice and the real flag manifold  as a generalization of the computation of cohomology of the flag manifold in the complex case.  In the complex case the cohomology ring of the flag manifold  can be described in terms  of polynomials and Weyl  group invariant polynomials on a Cartan subalgebra.   In the real split case, a split  Cartan  subgroup $H$ continues to play a central role but  now it is regarded as the \lq\lq  isospectral manifold of the Toda flow\rq\rq.   To be precise, there is a Toda flow on the connected component of the identity (\cite{kostant:79}), and also in (usually) disconnected open dense subsets of the remaining connected components of the Cartan (see  \cite{casian:02}).  This Toda flow depends,  in principle,  on a choice  $\gamma=(\gamma_1, \ldots, \gamma_l)\in{\mathbb R}^l$ of values for $l$  integrals of motion, i.e.
the Chevalley invariants, but its topological structure is independent of the value of $\gamma$ (if it is generic, i.e. ad-semisimple case with distinct real eigenvalues).   This Toda flow gives additional structure to the connected components of $H$, and this structure will play a role in the cohomology of the real flag manifold. The points in $H$ which are not in the isospectral manifold are blow-up points of the Toda flow. The set of blow-up points is called the Painlev\'e divisor
\cite{flaschka:91}, and it is defined by the zero set of the $\tau$-functions of the Toda lattice,
${\mathcal D}_j:=\{\tau_j(t_1,\ldots,t_l)=0\}$ for each $j=1,\ldots,l$ where $t_k$ are the flow parameters
of the Toda lattice (see below for the details). 
\smallskip

 We let  $\Pi=\{ \alpha_1,\ldots, \alpha_l \}$ denote simple roots relative to the Lie algebra ${\mathfrak h}$ of $H$ and we denote the root characters by $\chi_{\alpha_i}$. We can now  decompose $H$ into connected components $H_\epsilon$ where $\epsilon=(\epsilon_1\ldots \epsilon_l)$ with $\epsilon_i\in \{ \pm \}$ and $H_\epsilon:= \{ h\in H \,|\, {\rm sgn} (\chi_{\alpha_i}(h))=\epsilon_i,~i=1,\ldots,l\,\}$.  It is possible to arrange that there are exactly $2^l$ connected components in the split Cartan subgroup,  but  for now,  this is not necessary. Roughly speaking, the integral cohomology of $G/B$ emerges when one attempts to count the number of connected components   in the isospectral manifold   of the Toda lattice within one of the  connected components of the split Cartan subgroup  which we will denote $H_{-}$, with  $H_{-}:=H_{\epsilon}$ and $\epsilon=(- \ldots -)$. Other connected components of the Cartan subgroup lead to cohomology of $G/B$ with local coefficients.  From a different point of view, the connected components of the Cartan subgroup are the interior of polytopes, denoted by  $\Gamma_\epsilon$,  whose vertices correspond to the orbit of the Weyl group. To a first approximation, the graph of incidence numbers of the flag manifold
(with edges corresponding to non-zero (co) boundaries on Bruhat cells) is obtained by considering
a graph $\mathcal G$ where any two Weyl group elements satisfying $w\le w'$ (Bruhat order) with the lengths, $l(w')=l(w)+1$, are connected by an edge  if they are in the same connected component of the polytope  when blow-up points are removed. The precise definition of the graph $\mathcal G$ is in Definition \ref{graph}. We have the following theorem which is restated in a more general form in Theorem \ref{mainT} and then proved in Section \ref{mainTproof}:
\begin{Theorem}\label{mainT1}
The graph $\mathcal G$ defined in terms of the blow-up points of the Toda lattice is also the graph of
incidence numbers for the integral cohomology of the real flag manifold $\check{\mathcal B}=\check G/\check B$ in terms of the Bruhat cells.
\end{Theorem}

\smallskip
The connection with cohomology extends to etale cohomology over a field of positive characteristic. If,  in the context of the Toda lattice,  one attempts to count the number of blow-up points (singularities) along the trajectories of the  Toda flow, one ends up obtaining the same numbers that appear in the calculation of Frobenius eigenvalues in the context of the flag manifold (see the example below for the simplest case of this).  In more explicit terms, this corresponds to a surprising relation that we find  between the multiplicity $d$ of the singularity of the union of the Painlev\'e divisors, ${\mathcal D}_0=\cup_{j=1}^l{\mathcal D}_j$, at the point $p_o$ of intersection of all the divisors,
i.e. $\{p_o\}:=\cap_{j=1}^l{\mathcal D}_j$, and the polynomial giving the order of  a finite Chevalley group
$K({\mathbb F}_q)$.  We consider the polynomial  in $q$ defined by $|K({\mathbb F}_q)|$
 (e.g. $K=SO(n)$ for type $A$). This polynomial  has only two non-negative roots,  namely $q=1$ with multiplicity $g$ and $q=0$ with multiplicity $r$. For example, $|SO(3;{\mathbb F}_q)|=q(q^2-1)$, i.e. $g=1$ and $r=1$.  We define  $\tilde p(q):=q^{-r} |K({\mathbb F}_q)|$. Then we obtain that the multiplicity $d$ of this singularity  at $p_o$  is given by the formula $d={\rm deg}(\tilde p(q))$. This is verified for all type of classical semisimple Lie algebras and type $G_2$
 (Proposition \ref{direct}). This number  $d$  also agrees with the total number of blow-ups along the trajectory of the Toda flow (Proposition \ref{totalBnumber}).
 The polynomial  $\tilde p(q)$ is recovered from the Toda lattice  as a  polynomial  $p(q)$ defined by  an alternating sum of the number of blow-up points along trajectories  of the Toda lattice (Definition \ref{pq}). It is then shown in Theorem \ref{Kpq} that $p(q)=\tilde p(q)$. The number $d$, in the end, has a  simple Lie-theoretic description  as the dimension of any Borel subalgebra of  ${\rm Lie}(K({\mathbb C}))$.  This is due to Theorem 9.3.4 of \cite{carter}.  Then our second main result is summarized as follows (which follows from Theorem \ref{Kpq} proved in Section \ref{coh}):
\begin{Theorem} \label{mainT2}
The polynomial $\tilde p(q)=q^{-r}|\check K({\mathbb F}_q)|$ given in terms of the order of the finite Chevalley group $\check K({\mathbb F}_q)$ can be recovered in terms of the blow-ups of the Toda lattice, that is, we have $p(q)=\tilde p(q)$.
\end{Theorem}
This theorem corresponds to a calculation of the cup product structure in Theorem \ref{mainT1}
but over the  rationals. This is due to the fact that $\tilde p(q)$ contains all the information on the cup product structure of the compact group $\check K$.
\smallskip

\begin{figure}[t]
\includegraphics[width=3.7in]{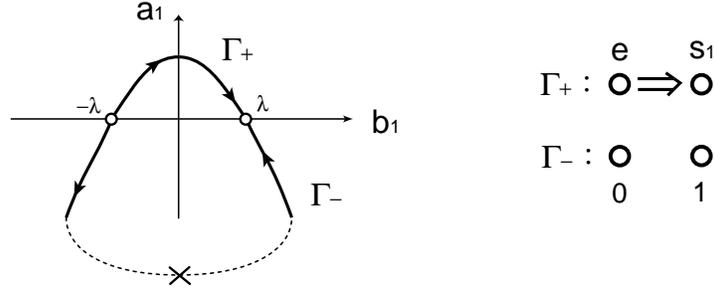}
\caption{The $\Gamma_{\epsilon}$-polytopes for the $A_1$-Toda lattice.
The left figure shows the invariant curve $I_1=a_1+b_1^2=\lambda^2$,
and the mark $\times$ indicates the blow-up point, the Painlev\'e divisor. The right figure shows the
graphs of $\Gamma_{\pm}$-polytopes with the Weyl action.
This shows that $\Gamma_+$ is connected, and $\Gamma_-$ has two connected components
separated by the Painlev\'e divisor.  The numbers in the graph of $\Gamma_-$
indicate $\eta(w)$, i.e. $\eta(e)=0$ and $\eta(s_1)=1$.}
\label{A1:fig}
\end{figure}

Let us illustrate our main results  by taking the simplest example ${\mathfrak g}={\mathfrak{sl}}(2;{\mathbb R})$, which also provides the basic structure of the general case. The ${\mathfrak{sl}}(2;{\mathbb R})$-Toda lattice is expressed by the Lax equation,
\[
\frac{dL}{dt}=[L,A]=LA-AL, 
\]
where the $2\times 2$ matrices $L$ and $A$ are given by
\[
L=\begin{pmatrix} b_1 & 1 \\ a_1 & -b_1\end{pmatrix}\,,\quad\quad
B=\begin{pmatrix} 0 & 0 \\ a_1 & 0 \end{pmatrix}\,.
\]
The flow of the Toda lattice is expressed by the invariant curve,
\[
I(a_1,b_1)=\frac{1}{2}{\rm Tr}(L^2)=a_1+b_1^2={\rm constant}.
\]
Then the isospectral manifold $Z(\gamma)_{\mathbb R}$ for a real split case is given by
the curve $I(a_1,b_1)=\gamma_1>0$ (see Figure \ref{A1:fig} where $\gamma_1=\lambda^2$). The compactified manifold
${\tilde{Z}}(\gamma)_{\mathbb R}$ consists of two polytopes (line segments) $\Gamma_{\epsilon}$
with $\epsilon={\rm sgn}(a_1)$, i.e.
\[
{\tilde Z}(\gamma)_{\mathbb R}=\Gamma_+\cup \Gamma_- \cong S^1\,,
\]
where $\Gamma_{\epsilon}=\overline{(\exp(tL^0)B/B)}$ with $L^0$ the initial
matrix of $L$ having ${\rm sgn}(a_1^0)=\epsilon$, and $B$ the Borel subgroup of upper
triangular matrices. The compactification is obtained by the companion embedding,
$Z(\gamma)_{\mathbb R}\to G/B$ (see \cite{flaschka:91, casian:02b}).
The $\Gamma_-$ has two connected components of the regular flows which are isolated 
by the singularity, called the Painlev\'e divisor marked by
$\times$ in the Figure \ref{A1:fig}.
We then define {\it graphs} associated with those polytopes $\Gamma_{\epsilon}$,
whose vertices are given by the fixed points, the end-points of $\Gamma_{\epsilon}$
and with an edge $\Rightarrow$ when the vertices
are connected by the flow, i.e. no blow-up (see Definition \ref{graph} for more details).
The $\Gamma_+$ gives a connected graph with one edge $\Rightarrow$, and
the graph of $\Gamma_-$ consists of two isolated vertices. 
We also assign the end-points (vertices) of $\Gamma_{\epsilon}$ with the elements $e$ and $s_1$
of the Weyl group $W=S_2$, the symmetry group of order two, that is, the vertices are given by
the orbit of the Weyl group (recall that the Toda flow defines a torus action on the flag manifold, and the compactified manifold $\tilde Z (\gamma)_{\mathbb R}$ is a smooth toric
variety \cite{casian:02}).
In particular, the graph of $\Gamma_-$ is the same as the incidence graph for the real
flag manifold $SL(2;{\mathbb R})/B\cong S^1$, where the edge has
the incidence number 2 originally proved in \cite{casian99}
(this is stated as Theorem \ref{mainT} which uses Theorem \ref{casian99theorem}). 

We now explain the connection between the blow-up of the Toda lattice and the cohomology
of a compact Lie group ($SO(2)\cong S^1$ in this case).
Let us first  consider  the complexification ${\mathbb C}^*={\mathbb C}\setminus \{ 0 \}$ of this flag manifold; and then  its reduction to positive characteristic $k_q^*=k_q\setminus \{ 0 \}$ where $k_q$ is an algebraic closure of the finite field ${\mathbb F}_q$ and $q$ is a power of a prime. Consider then the set of ${\mathbb F}_q$ points on this flag manifold. This is the number of elements in ${\mathbb F}_q\setminus \{ 0 \}$, namely $q-1$. Similarly, we can just consider the number of elements  $|S^1({\mathbb F}_q)|$, that is, pairs $(x,y) \in {\mathbb F}_q\times {\mathbb F}_q$ satisfying the equation $x^2+y^2=1$. The number $|S^1({\mathbb F}_q)|$  gives a polynomial in $q$ which can also be obtained cohomologically by considering  
the Frobenius  map, $Fr_q: k_q^* \to k_q^*$, $z\mapsto z^q$ whose Lefschetz number
(alternating sum of the traces of $Fr_q$ in (etale) cohomology with proper supports) is given by
\[
L(Fr_q)=q-1.
\]

A characteristic zero analogue of this positive characteristic construction  can be given in this simple case.  We consider the map:
$\Phi_q:S^1 \to S^1$, $z\mapsto z^q$, a map of degree $q$.  Then we have
\[
\left\{\begin{array}{llll}
H_0(S^1;{\mathbb Q})={\mathbb Q}\quad &:\quad q^0\\ 
H_1( S^1;{\mathbb Q})={\mathbb Q}\quad &:\quad q^1
\end{array}\right.
\]
and now the number of fixed points is $q-1$, i.e. the number of non-zero roots of $z^q=z$, or as in
the Lefschetz fixed point theorem,
\[
L(\Phi_q)={\rm Tr}\left((\Phi_q)_*|_{ H_1( S^1;{\mathbb Q})}\right)- {\rm Tr}\left((\Phi_q)_*|_{H_0(S^1;{\mathbb Q})}\right)=q-1\,.
\]
 We now  obtain this  polynomial, given as a  Lefschetz number, in terms of the singularities of the  $A_1$   Toda flow.  Assign  to each vertex $w$ in the graph of $\Gamma_-$  the number of singularities (blow-ups) counted along the flow starting from the top vertex $e$ to 
the vertex $w$: We denote the number by $\eta(w)$, so that $\eta(e)=0$. This polynomial is introduced in Definition \ref{pq} as
\[
p(q):=-\sum_{w\in S_2}(-1)^{l(w)}q^{\eta(w)}=-( (-1)q^{\eta(s_1)}+q^{\eta(e)})=q-1,
\]
which agrees with $L(Fr_q)$, that is, we have (Theorem \ref{Kpq})
\[
p(q)=L(Fr_q)\,.
\]
If $q$ is a power of an odd prime such that $x^2+1$ factors over  ${\mathbb F}_q$,  this polynomial also gives the number of ${\mathbb F}_q$ points
on $SO(2) \cong S^1$, i.e.
\[
|S^1({\mathbb F}_q)|=q-1.
\]
This can be directly computed from the stereographic projection of the circle. Namely  with
$x=2s/(s^2+1),~y=(s^2-1)/(s^2+1)$, then count the number of solutions $s\in {\mathbb F}_q$ for
$x^2+y^2=1$. If $s^2+1=0$ factors over ${\mathbb F}_q$ then there are $q-2$ values of $s$ that give rise to $x,y$ in these formulas plus the point $(0,1)$. This gives a total of  $q-1$ points.
Also note that the compact group $S^1\cong SO(2)$ is the maximal compact
subgroup of $G=SL(2;{\mathbb R})$ (Theorem \ref{Kpq}). The main purpose of the present paper is to
clarify the correspondence among those objects for the general case of $\mathfrak g$-Toda lattice with real split
semisimple Lie algebra ${\mathfrak g}$.

We also mention some analytic structure of the Toda lattice in the connection of
the blow-up with the degree of polynomial $p(q)$:
The solution $(a_1,b_1)$ can be expressed by the $\tau$-function (see (\ref{ab}) and (\ref{tau}) for
the details),
\[
a_1(t)=\frac{\pm 1}{\tau_1(t)^2},\quad b_1(t)=\frac{d}{dt}\ln\tau_1(t)\,.
\]
The $\tau$-function, $\tau_1(t)$, is given as follows:
For the case $a_1>0$ (i.e. $\Gamma_+$-polytope), the function $\tau_1$ can be calculated
by the companion matrix $C_{\gamma}$ of $L$ with $\gamma=\lambda^2$ for $a_1^0=1$,
\[
\tau_1(t_1)=\langle \exp(tC_{\gamma}) e_1,e_1\rangle={\rm cosh}(\lambda t_1),
\]
where $C_{\gamma}=\begin{pmatrix}0&1\\\lambda^2&0\end{pmatrix}$
and $e_1=(0,1)^T$ (see (\ref{tau}) and also \cite{casian:04}).
Since there is no blow-up in this case, there is an edge between the vertices in the graph  of $\Gamma_+$ (see Figure \ref{A1:fig}).
For the case $a_1<0$ (i.e. $\Gamma_-$-polytope), we have
\[
\tau_1(t_1)=\langle\exp(tC_{\gamma})e_2,e_1\rangle=\frac{1}{\lambda}{\rm sinh}(\lambda t_1),
\]
where $e_2=(0,1)^T$.
Thus the solution $(a_1(t), b_1(t))$  blows up at $t_1=0$. Near $t_1=0$, the $\tau_1$-function has the form, $\tau_1=t_1+.....$,
in which the first term is the Schur polynomial $S_{(1)}(t_1)$ (see (\ref{schur-tau}) for the $A_l$-Toda lattice). 
The number of zeros in the $\tau$-functions gives the total number of blow-ups in the Toda flow
from the top vertex to the bottom one. 
 In the case of ${\mathfrak{sl}}(2;{\mathbb R})$-Toda lattice, this number equals to one which agrees with the degree of polynomial $p(q)=q-1$. The agreement that the total number of blow-ups in the Toda
 flow is equal to the degree of the polynomial $p(q)$ for corresponding Chevalley group $K({\mathbb F}_q)$ is true for the general case.
 In Proposition \ref{totalBnumber}, we give the explicit numbers for the general case. 
 We also show in Proposition \ref{direct} that in all types of $\mathfrak g$ except $E$ and $F$ the degree of the polynomial $p(q)$ is given by the sum of the minimum degrees in the $\tau$-functions near the point $\{p_o\}=\cap_{j=1}^l{\mathcal D_j}$.
 We then conjecture that the degree is also equal to the number of {\it real} roots of the product of $\tau$-functions (Conjecture \ref{tangentcone}). The number of the complex
roots has been found to be the height $|2\rho|$ where $\rho$ is the sum of all the fundamental weights
(\cite{adler:91, flaschka:91}).

The numbers $\eta(w)$ of blow-ups along the Toda flow are introduced in Definition \ref{eta}, and the maximum number 
is obtained by $\eta(w_*)$ with the longest element $w_*$. Our third main result is the following
(Propositions \ref{direct} and \ref{totalBnumber}):
\begin{Theorem}
We have $\eta(w_*)={\rm deg}(\tilde p(q))$, the dimemsion of any Borel subalgebra of ${\rm Lie}(\check K({\mathbb C}))$.
Assuming that the Lie algebra ${\mathfrak g}$ has no factor of type $E$ or $F$, then we have
$\eta(w_*)=d$, where $d$ is the multiplicity of the sigularity of ${\mathcal D}_0=\cup_{j=1}^l{\mathcal D}_j$ at the point $p_o$ of the intersection of all the Painlev\'e divisors ${\mathcal D}_j=\{\tau_j(t_1,\ldots,t_l)=0\}$.
\end{Theorem}

\smallskip

Finally we remark that in the simple example of $SL(2;{\mathbb R})$, the connection with cohomology and representation theory of the group can be seen directly.  We consider the de Rham chain complex on $G/B\cong K/T\cong S^1$. Let $F(S^1)$ and $\Omega(S^1)$ be the ${\mathbb C}$-valued functions and $1$-forms on $G/B$ with a finite Fourier expansion on $S^1$. The cohomology of $S^1$ appears from the chain complex  $M:=F^{even}(S^1)\to \Omega^{even}(S^1)=:\tilde M$   involving only even powers of $e^{\sqrt{-1}\theta}$.  Following the general arguments in \cite{casian99}, this is a map of ${\mathfrak{sl}}(2;{\mathbb R})$-modules and the two modules are principal series representations. Note  that  $F^{even}(S^1)$ consists of
functions on $G/B\cong K/T$ and contains constant functions as a trivial  ${\mathfrak g}$-submodule and $\Omega^{even}(S^1)$
which consists of $1$-forms only contains the trivial ${\mathfrak g}$-module as a quotient. 
The structure of these modules can be described explicitly as well as the map in the chain complex.  
 There are three irreducible modules $C$, $D_+$, $D_-$. The module $C$ is  one dimensional and corresponds to the constant functions; each of $D_\pm$ is a discrete series representation.  The modules $D_{\pm}$ consist of functions involving linear combinations of $e^{ n\theta \sqrt {-1}}$ with ${\rm sgn}(n)=\pm $. We use  labels  $q^{i/2}$ to distinguish submodules from quotients, on which  the map of ${\mathfrak g}$-modules (coboundary)
sends the quotient $D_+\oplus D_-$ of $M$ to the submodule $D_+\oplus D_-$ of $\tilde M$:
$$ \begin{matrix} 
q^{1/2}&:&                    &{D_+\oplus D_-}& {}               & {} & {C}& :& q^1 \\
{}         & &            &{} & {\searrow }                  & {} &{}&   \\
 q^0 &:&             &{C}& {}                  & {} &{D_+\oplus D_-}& :& q^{1/2} \\
 \end{matrix} $$
Note here that we have normalized the labels in $\tilde M$, so that $D_+\oplus D_-$,
which is in the image of the map,  carries the same label in both modules. This map has $q^0C$
in the kernel and $q^1C$ as cokernel giving the cohomology of the circle: 
$H^*(S^1;{\mathbb C})=\Lambda(q^1x_1)$, i.e.  $H^0(S^1;{\mathbb C})=C$ (with label $q^0$)
and $H^1(S^1;{\mathbb C})=C$ (with label $q^1$). This corresponds to the polytope $\Gamma_{-}$ and the graph with no edges between $e$ and $s_1$ (Theorem \ref{mainT}).

  The filtration described
 (e.g. $C \subset M$ and $D_+\oplus D_-\subset \tilde M$) is an example of a
 {\it weight filtration} and arises from a filtration via Frobenius eigenvalues (see \cite {casian86} for
a description of these filtration of principal series modules).
The labels $q^i$ attached to the cohomology are Frobenius eigenvalues which appear when one considers a field of positive characteristic  $k$, then replace
$k\setminus \{ 0 \}$ instead of ${\mathbb C}^*$,  and use  etale cohomology with coefficients
$\bar {\mathbb Q}_m$ for appropriate $q$ and $m$.  The relation to the order of the finite group $|SO(2;{\mathbb F}_q)|=q-1 $ follows by considering the alternating sums of  the trace of Frobenius  in etale cohomology with proper supports.

There is another chain complex $F^{odd}(S^1)\to \Omega^{odd}(S^1)$ which computes cohomology  of $G/B\cong S^1$ with twisted coefficients.  The coboundary here is an isomorphism,  and  therefore the rational cohomology with twisted coefficients is zero.  This corresponds to the $\Gamma_+$ polytope and the graph $e\Rightarrow s_1$. The edge $\Rightarrow$ represents a coboundary and $H^*(S^1;{\mathcal L} )=0$ for appropriate twisted coefficients  ${\mathcal L}$.


\section{Preliminaries and notations}
We now introduce some of the main objects needed to state our main resuls.
  We will use standard Lie group notation: $B, B^*$ are
a pair of opposite Borel subgroups in a split semisimple Lie group $G$,   $B=HN$, $B^*=HN^*$, $H$ a split Cartan subgroup,
$W$ the Weyl group. The Lie algebras corresponding to those are denoted as
${\mathfrak b}={\mathfrak h}\oplus{\mathfrak n}$ and ${\mathfrak b}^*={\mathfrak h}\oplus
{\mathfrak n}^*$. The group $G$ is more carefully defined below.

\subsection{The real flag manifold and the group $G$}
\label{group2}
 Let ${\mathfrak g}$ be real split semisimple Lie algebra, $ {\mathfrak g}_{\mathbb C}$ its complexification. We now consider any complex connected Lie group $G({\mathbb C})$ with Lie algebra ${\mathfrak g}_{\mathbb C}$.  The real  Lie subalgebra  ${\mathfrak g}$ corresponds uniquely  to a connected Lie subgroup  $G$.  
 
 We can now regard $G({\mathbb C})$ as an algebraic group and consider the group of real points $G({\mathbb R})$ containing $G$ as a connected component.  We  have an
algebraic group defined over a finite extension of  ${\mathbb Q}$. 
For any field $k$, we  denote by  $G(k)$  the set of  $k$ points.
This immediately makes sense when $k={\mathbb C}$ or ${\mathbb R}$,
but as is pointed out below,  by reducing to positive characteristic,  we can also make sense of $G(k)$ when $k$  has characteristic of a prime $p$.  

We  recall that  $\check {\mathfrak g}$ the Lie algebra whose Cartan matrix is the transpose of the Lie algebra  $ {\mathfrak g}$. The previous construction then will give rise to a connected group $\check G$ contained in the simply connected complex group $\check G({\mathbb C})$. We also have flag manifolds ${\mathcal B}=G/B$ and $\check {\mathcal B}={\check G}/{\check B}$. 

\begin{Remark}\label{group1}  Let $G({\mathbb C})$ be simply connected. Then the  group  ${\rm Ad}(G)$ is what is denoted $G$ in \cite{casian:02}.   The group $\tilde G=\{ {\rm Ad}(g) \in {\rm Ad}(G({\mathbb C}))\,|\,{\rm Ad}(g){\mathfrak g} \subset {\mathfrak g} \} $  is introduced in \cite{casian:02} because its  split
  Cartan subgroup $H_{\mathbb R}$   with Lie algebra ${\mathfrak h}$ has exactly $2^l$ connected
components (matching the $2^l$ polytopes arising from the Toda lattice).  Then, using the work of Kostant in \cite{kostant:79},  we showed in \cite{casian:02}
that in the regular  ad-semisimple case  the isospectral manifold of the Toda lattice is embedded as an open dense set in $H_{\mathbb R}$.  
The connected components of $H_{\mathbb R}$ can  be pictured as the interior 
of polytopes $\Gamma_\epsilon$ with vertices corresponding to the Weyl group $W$ of ${\mathfrak g}$. 
For example the four hexagons in the case
of $A_2$ (with blow-ups included) are identifiable with the connected components
of the Cartan subgroup in $SL(3;{\mathbb R})$ consisting of diagonal matrices (see also Figure \ref{hexagon1:fig}).
\end{Remark}

We let $K$ denote a maximal compact subgroup of $G$. It is convenient
to think of $K$ as the fixed point set of an involution $\theta$ (Cartan involution). Hence
we consider an algebraic group  endowed with an involution $\theta$ and $K(k)$ consists of the $\theta$ fixed points of $G(k)$. 

\begin{Example} We consider $G({\mathbb C})=SL(n;{\mathbb C})=\check G ({\mathbb C})$, i.e. the set 
of all the $n\times n$ complex matrices $A$ satisfying
the polynomial equation, ${\rm det}(A)=1$.  Then the complex
solutions of  ${\rm det}(A)=1$ define $G({\mathbb C})$. The group
of real points is, of course, $G({\mathbb R})=SL(n;{\mathbb R})$.
The involution $\theta$ is given by $\theta(A)=A^*$, the inverse of the
transpose of $A$. We then have $K=SO(n)$ as the set of matrices satisfying $\theta(A)=A$. 
\end{Example}

\subsection{Reduction to positive characteristic}

We will show how   the Toda lattice and the structure of its blow-up points  contains the cohomology of the real flag manifold and of the compact group $K$.  It turns out that it  also contains some additional structure which is  present  when the flag manifold is regarded over a field of positive characteristic.   
We then need to make sense of the main Lie-theoretic objects over a field ${\mathbb F}_q$ with $q$ elements. 

\begin{Notation}
 We denote $k_q=\bar{\mathbb F}_q$
an algebraic closure of the field ${\mathbb F}_q$ with $q$ elements.
\end{Notation}

The algebraic groups considered above  as well as the involutions $\theta$ are
defined over a finite extension of ${\mathbb Q}$. Therefore
it is possible to reduce to positive characteristic and define 
 an algebraic group defined over a field ${\mathbb F}_q$
with $q$ elements, where $q$ is a power of a prime number $p$, $p\not=2$.
We then have the group $G(k_q)$.  The involution
$\theta$ gives rise to an involution $\theta$ defined over ${\mathbb F}_q$ and this 
leads to $K(k_q)$, the points in 
$G(k_q)$ fixed by the involution $\theta$. We also obtain
the finite Chevalley groups $G({\mathbb F}_q)$ and $K({\mathbb F}_q)$
by considering the subsets of ${\mathbb F}_q$ points.   These two finite groups are, respectively,  the fixed points of the Frobenius map   $Fr:G(k_q)\to G(k_q)$  and  $Fr:K(k_q)\to K(k_q)$ which are reviewed in the example below.

\begin{Example} \label{det}
The equation ${\rm det}(A)=1$ has integral coefficients.  The integral coefficients
 make it possible to reduce modulo a prime $p$. Let $q$ be a power of a prime
$p$ and $k_q$ an algebraic closure of ${\mathbb F}_q$. We may then consider  $G(k_q)=SL(n;k_q)$ the solutions in $k_q$ of ${\rm det}(A)=1$.   The Frobenius map is given by $Fr( (a_{ij}))=( a_{ij}^q)$.
By reducing modulo $p$ the involution $\theta(A)=A^*$
we obtain $SO(n;k_q)$ as the set of 
fixed points and then the finite group  $SO(n;{\mathbb F}_q)$.
In the simplest example, $n=2$ and  $SO(n;{\mathbb F}_q)$  consists of $2\times 2$ matrices
$\left(
\begin{matrix}
x & y   \\
-y & x \\
\end{matrix}\right)$ satisfying $x^2+y^2=1$ with $x,y \in k_q$.
\end{Example}

\begin{Remark} \label{so}
The number of points of $SO(n;{\mathbb F}_q)$ can be found in  p.75 of \cite{carter}.  The formula depends on the  Witt index,  the dimension of a maximal isotropic subspace. For instance over $\mathbb F_3$ the Witt index of $SO(2;{\mathbb F}_3)$ is zero  and $|SO(2;\mathbb F_3)|=q+1=4$;  however over characteristic $5$ there is an isotropic vector space since $1^2+2^2=0$. In this case  $|SO(2;\mathbb F_5)|=q-1=4$.  By extending $\mathbb F_3$ to  $\mathbb F_{3^2}$ (adding $\sqrt{-1}$)   the  number of ${\mathbb F}_q$ points becomes $q-1$.   In general we have  
$|SO(n+1;{\mathbb F}_q)|=  |SO(n;{\mathbb F}_q)|\times|S^{n}({\mathbb F}_q)|$ (with $S^n$ the sphere given by $x_1^2+\cdots +x_{n+1}^2=1$). Hence
$|SO(n+1;{\mathbb F}_q)|=|S^{1}({\mathbb F}_q)|\times |S^{2}({\mathbb F}_q)| \times \cdots  \times |S^{n}({\mathbb F}_q)|$. For example   $|SO(6;{\mathbb F}_q)|= |SO(5;{\mathbb F}_q)|\times |S^{5}({\mathbb F}_q)|$.  We have $|SO(5;{\mathbb F}_q)|=q^4(q^2-1)(q^4-1)$. However we find that for $q=3$, $|S^{5}({\mathbb F}_q)|=252=q^2(q^3+1)$ and for $q=5$, $|S^{5}({\mathbb F}_q)|=3100=q^2(q^3-1)$. 
This gives two different formulas for  $|SO(6;{\mathbb F}_q)|$. This is because when $q=3$ one obtains a group of type ${}^2A_3$ while $q=5$ gives a group of type $A_3$ (see \cite{carter}). The problem is the same as in the case of $SO(2;\mathbb F_q)$ and it disappears  by considering $q^2$ for $q=3$.

We illustrate the computation of $|S^{n}({\mathbb F}_q)|$ when $\sqrt{-1} \in {\mathbb F}_q$ by doing the case of $n=2$. In particular we compute $|SO(3;{\mathbb F}_q)|$:
 Using the formulas for the  stereographic projection;  $x=\frac{2u}{1+u^2+v^2}$, 
 $y=\frac{2v}{1+u^2+v^2}$, $z=\frac{u^2+v^2-1}{1+u^2+v^2}$ and 
 $u=\frac{x}{1-z}$, $v=\frac{y}{1-z}$, we have a bijective correspondence between the set  $\{ (x,y,z) \in {\mathbb F}_q^3 \,|\, z\not= 1 \}$ and $\{ (u,v)\,|\, 1+u^2+v^2\not =0 \}$.
 Multiplying by  $\sqrt {-1}\in {\mathbb F}_q$, this second set also corresponds bijectively to $\{ (u,v)\,|\,u^2+v^2=1 \}$.
 This produces $q^2-(q-1)$ points.  We now must compute the \lq\lq north pole\rq\rq points, i.e. 
 all the $(x,y,1)$ in the sphere. We then have $x^2+y^2=0$ and thus for $y\not=0$, $x=\pm \sqrt{-1}y$, that is, $2(q-1)$ points. Adding $(0,0,1)$ we obtain a total of $2(q-1)+1=2q-1$ \lq\lq north pole\rq\rq elements.  We obtain a total of $q^2-(q-1)+ 2q-1=q^2+q=q(q+1)$.
Thus, for example, $SO(3;{\mathbb F}_q)$ has 
 $|S^1({\mathbb F}_q)|\times  |S^2({\mathbb F}_q)|=(q-1)q(q+1)=q(q^2-1)$ points. 
 \end{Remark}


 We now introduce the \lq\lq complex version of the real flag manifold\rq\rq.  Consider the $K({\mathbb C})$-orbit ${\mathcal O}_o$ which is the unique
 open dense orbit in  the complex flag manifold ${\mathcal B}_{\mathbb C}=G_{\mathbb C}/B_{\mathbb C}$. This orbit has the
 homotopy type of ${\mathcal B}=G/B$ and can be regarded as a {\it complex
 version} of the real flag manifold ${\mathcal B}$. This is the object which we actually
study, not just over ${\mathbb C}$ but  also over  a field  of positive characteristic.  We also consider all the   $K({\mathbb C})$-equivariant local systems ${\mathcal L}$  over  ${\mathcal O}_o$. These local systems contain 
 the information to compute integral cohomology of the real 
 flag manifold (\cite{casian99}).  
    
 
\begin{Example}\label{suoneone} Let  $G(k)= SL(2; k)$.  We consider
the Cartan involution  $\theta$ given by
 $\theta(g)= \left(
\begin{matrix}
1 & 0   \\
0 & -1 \\
\end{matrix}\right) g \left(
\begin{matrix}
1 & 0   \\
0 & -1 \\
\end{matrix}\right)$. Therefore $K(k)$ as the set of fixed points of $\theta(g)=g$ consists of the diagonal matrices  $\left\{ {\rm diag}(z, z^{-1}) \,| \,z\in k^* \right\}$. The choice of $K$ corresponds to considering
 the subgroup $SU(1,1)$ of $SL(2; {\mathbb C})$ rather than
$SL(2; {\mathbb R})$.  The flag manifold is ${\mathbb P}^1 = k \cup \{ \infty \}$ and the action of 
$K(k)$ is $g(z)\cdot y = z^2y$ for $g(z)\in K(k)$ with $z\in k^*$. The $K(k)$-orbits are $k^*$, $\{ 0 \}$ and $\{ \infty \}$.
Hence for $k={\mathbb C}$, we have ${\mathcal O }_o={\mathbb C}^*$. This orbit has the homotopy type of 
the real flag manifold of $SU(1,1)$, namely a circle $S^1$. Take $k=k_q$ and note the  map given by $Fr(z)=z^q$. The fixed points are the ${\mathbb F}_q$ points.  The map $Fr$ then induces multiplication by $q$ on $H^1( k_q\setminus \{ 0 \} ;\bar {\mathbb Q}_m)$ (etale cohomology). 

\end{Example}

\subsection {The Toda lattice}
Here we present
the background information on the Toda lattice:
The generalized (nonperiodic) Toda lattice equation related to
real split semisimple Lie algebra $\mathfrak g$ of rank $l$ is defined by
the Lax equation,
(\cite{bogoyavlensky:76,kostant:79}),
\begin{equation}
\label{lax}
\displaystyle{\frac{dL}{dt}=[L,A]}
\end{equation}
where $L$ is a Jacobi element of ${\mathfrak g}$ and $A$ is the ${\mathfrak n}^*$-projection
of $L$, denoted by $\Pi_{\mathfrak n^*}L$,
\begin{equation}\left\{
\begin{array} {ll}
& \displaystyle{L(t)=\sum_{i=1}^l b_i(t)h_{\alpha_i}+\sum_{i=1}^l \left(
a_i(t)e_{-\alpha_i}+e_{\alpha_i}\right)} \\
& \displaystyle{A(t)=\Pi_{\mathfrak n^*}L=\sum_{i=1}^l a_i(t)e_{-\alpha_i}}
\end{array}
\right.
\label{LA}
\end{equation}
Here $\{h_{\alpha_i},e_{\pm\alpha_i}\}$ is the Cartan-Chevalley basis of $\mathfrak g$ with the simple roots
$\Pi=\{\alpha_1,\cdots,\alpha_l\}$
which satisfy the relations,
\[
    [h_{\alpha_i} , h_{\alpha_j}] = 0, \quad
    [h_{\alpha_i}, e_{\pm \alpha_j}] = \pm C_{j,i}e_{\pm \alpha_j} \ , \quad
    [e_{\alpha_i} , e_{-\alpha_j}] = \delta_{i,j}h_{\alpha_j},
\]
where $(C_{i,j})$ is the $l\times l$ Cartan matrix of
$\mathfrak g$.
The Lax equation (\ref{lax}) then gives
\begin{equation}
\left\{
\begin{array} {ll}
& \displaystyle{\frac{db_i}{dt}=a_i\,,} \\
& \displaystyle{\frac{da_i}{dt}=  -\left(\sum_{j=1}^l C_{i,j} b_j\right)a_i\,.}
\end{array}
\right.
\label{toda-lax}
\end{equation}

The integrability of the system can be shown by the existence of the Chevalley
invariants, $\{I_k(L)\,| \,k=1,\cdots,l\}$, which are given by the weighted homogeneous
polynomial
of $\{(a_i, b_i)\, |\, i=1,\cdots,l\}$ with ${\rm deg}(a_i)=2$ and ${\rm deg}(b_i)=1$.
Those invariant polynomials also define the commutative equations
of the Toda equation (\ref{lax}),
\begin{equation}
\label{higherflows}
\frac{\partial L}{\partial t_k}=[L,\Pi_{\mathfrak n^*}\nabla I_k(L)]\,
\quad {\rm for}\quad k=1,\cdots,l\,,
\end{equation}
where $\nabla$ is the gradient with respect to the Killing form, i.e.
for any $x\in {\mathfrak g}$, $dI_k(L)(x)=K(\nabla I_k(L),x)$.
For example, in the case of ${\mathfrak g}={\mathfrak{sl}}(l+1,{\mathbb R})$,
the invariants $I_k(L)$ with ${\rm deg}(I_k)=k+1$ and the gradients $\nabla I_k(L)$ are given by
\[I_k(L)=\frac{1}{k+1}{\rm Tr}(L^{k+1})\,\quad{\rm and}\quad
\nabla I_k(L)=L^k.\]
The set of commutative equations is called the Toda lattice hierarchy.
Note that (\ref{toda-lax}) is the first member of the hierarchy, i.e. $t=t_1$.
Then the {\it real} isospectral
manifold is defined by
\[
Z(\gamma)_{\mathbb R}=\left\{(a_1,\cdots,a_l,b_1,\cdots,b_l)\in{\mathbb
R}^{2l}~\Big|~
I_k(L)=\gamma_k\in {\mathbb R},~ k=1,\cdots,l\right\}=\bigcap_{k=1}^l I_k^{-1}(\gamma_k)\,.
\]
The manifold $Z(\gamma)_{\mathbb R}$ can be compactified by adding the set of
points corresponding to the singularities ({\it blow-ups}) of the solution.
The compactification is obtained through the companion embedding, $Z(\gamma)_{\mathbb R}\to G/B$ (see \cite{flaschka:91,casian:02b}).
Then the compact manifold ${\tilde Z}(\gamma)_{\mathbb R}$ is
described by a union of convex
polytopes $\Gamma_{\epsilon}$ with $\epsilon=(\epsilon_1,\cdots,\epsilon_l),~
\epsilon_i={\rm sgn}(a_i)$ \cite{casian:02},
\[
{\tilde Z}(\gamma)_{\mathbb R}=\bigcup_{\epsilon\in\{\pm\}^l}\Gamma_{\epsilon}.
\] 
Each polytope
$\Gamma_{\epsilon}$ is expressed as the closure of the orbit of a
Cartan subgroup in the flag manifold, i.e.
\[
\Gamma_{\epsilon}=\overline{(G^{L^0}B/B)},
\quad{\rm with} \quad \epsilon=(\epsilon_1\ldots\epsilon_l), ~\epsilon_i={\rm sgn}(a_i^0)\,,
\]
where $L^0$ is an initial matrix of $L(t)$, and the Cartan subgroup is given by the orbit,
\[
G^{L^0}:=\left\{~\exp\left(\sum_{k=1}^l t_k\nabla I_k(L^0)\right)~\Big|~t_k\in{\mathbb R},~k=1,\ldots,l\right\}\,.
\]
In particular, the polytope $\Gamma_{-\ldots-}$, denoted simply as $\Gamma_-$, is given by 
\begin{equation}
\label{gcorbit}
\Gamma_{-}=\overline{(G^{C_{\gamma}}w_*B/B)}\,,
\end{equation}
where $C_{\gamma}$ is the companion matrix of $L$ and $w_*$ is the longest element
of the Weyl group $W$.  All $\Gamma_{\epsilon}$ are the same
polytope with the vertices given by the orbit of Weyl group (see \cite{atiyah:82}).
Each $\Gamma_{\epsilon}$ has a unique intersection with the Painlev\'e divisors
(see below and \cite{casian:02b}).
Thus in an ad-semisimple
case with distinct real eigenvalues, the
compact manifold ${\tilde Z}(\gamma)_{\mathbb R}$ is a toric variety, and
the vertices of each polytope are labeled by the elements of the Weyl group
(see Figures below for the examples of $\Gamma_{\epsilon}$ for types
$A_1, A_2$ and $G_2$).

\subsubsection{The $\tau$-functions and Painlev\'e divisors}
The solution of the Toda lattice is obtained by the
$\tau$-functions,
which are defined by
\begin{equation}
\label{ab}
b_k=\displaystyle{{\frac{d}{dt}}\ln \tau_k, \quad \quad a_k=a_k^0\prod_{j=1}^l
(\tau_j)^{-C_{k,j}}},
\end{equation}
where $a_k^0$ are some constants (those are obtained from (\ref{toda-lax})).
  The $\tau$-functions are given by \cite{flaschka:91}
\begin{equation}
\label{tau}
\tau_j(t_1,\ldots,t_l)=\langle g(t_1,\ldots,t_l)\cdot v^{\omega_j},v^{\omega_j}\rangle,\quad
g\in G^{L^0}\,.
\end{equation}
Here $v^{\omega_j}$ is the highest weight vector in the fundamental representation
of $G$, and $\langle\cdot,\cdot\rangle$ is a pairing on the representation space.
The blow-up points are given by the zeros of the $\tau$-functions, $\tau_j(t_1,\ldots,t_l)=0$
for some $j\in\{1,\ldots,l\}$. We then define the Painlev\'e divisors as the sets of zeros of $\tau$-functions:
\begin{Definition} (\cite{casian:02})
 The Painlev\'e divisor ${\mathcal D}_J$ for
a subset $J\subset\{1,\ldots,l\}$ is defined by
\[
{\mathcal D}_J:=\bigcap_{j\in J}{\mathcal D}_j\quad {\rm with}\quad {\mathcal D}_j:=\{\,(t_1,\ldots,t_l)\,|\,\tau_j(t_1,\ldots,t_l)=0\,\}\,.
\]
\end{Definition}
The ${\mathcal D}_J$ can be also described by the intersection with the Bruhat cell
$N^*w_JB/B$ with the longest element $w_J$ of the Weyl subgroup $W_J=\langle
s_j \,|\, j\in J\rangle$ (\cite{flaschka:91,adler:91}).
In particular, the divisor ${\mathcal D}_{\{1,\ldots,l\}}$ is a unique point, denoted as $p_o$, in the variety
$\tilde Z(\gamma)_{\mathbb R}$, and it is contained in the $\Gamma_{-}$-polytope (see \cite{casian:04}). 

In the case of $A_l$-Toda lattice, one can express the geometric structure of the
divisor ${\mathcal D}_0:=\cup_{j=1}^l{\mathcal D}_j$ as an algebraic variety near this point $p_o$
 (other cases including $B_l$ and $C_l$ are also discussed in \cite{casian:04}):
Setting this point as $t=(0,\ldots,0)\in {\mathbb R}^l$, the $\tau$-functions
are given by the following expansion around this point (use (\ref{tau}) with $v^{\omega_j}=e_1\wedge\cdots\wedge e_j$ and $\{e_k\}$ the standard basis of ${\mathbb R}^l$),
\begin{equation}
\label{schur-tau}
\tau_k=(-1)^{\frac{k(k-1)}{2}}S_{(l-k+1,\ldots,l)}(t_1,\ldots,t_l)+ ({\rm higher~degree~terms})\,.
\end{equation}
where $S_{(i_1,\ldots,i_k)}$ with $1\le i_1<\cdots<i_k\le l$ is the Schur polynomial defined as
the Wronskian determinant with respect to the $t_1$-variable,
\[
S_{(i_1,\ldots,i_k)}={\rm Wr}(h_{i_1}, h_{i_2},\ldots,h_{i_k})=\Big|\left(h_{i_{\alpha}-\beta+1}\right)_{1\le\alpha,\beta,\le k}\Big|\,.
\]
Here $h_j(t_1,\ldots,t_l)$ are complete homogeneous symmetric functions defined by
\[
\exp\left(\sum_{k=1}^{l} \lambda^k t_k\right)=\sum_{k=0}^{\infty}\lambda^k h_k(t_1,\ldots,t_l)\,
\quad {\rm with}\quad h_n=0~{\rm when}~n<0.
\]
Explicitly we have (see \cite{macdonald:79}),
\begin{equation}
\label{schurp}
\begin{array}{lllll}
h_k(t_1,\cdots,t_k)&=&\displaystyle{\sum_{k_1+2k_2+\cdots + nk_n=k} 
\frac{t_1^{k_1}t_2^{k_2}\cdots t_n^{k_n}}{k_1!k_2!\cdots k_n!}}\\
&{}&\\
&=&\displaystyle{\frac{t_1^{k}}{k!}+\frac{t_1^{k-2}t_2}{(k-2)!}+\cdots+t_{k-1}t_1+t_k\,.}
\end{array}
\end{equation}
In particular, $S_{(k)}=h_k$.
Those Schur polynomials $S_{(l-k+1,\ldots,l)}$ are the $\tau$-functions for the nilpotent Toda lattice for $A$-type \cite{casian:04},
and the Young diagram associated with $S_{(l-k+1,\ldots,l)}$ is the rectangular box
with $k$-stacks of $(l-k+1)$-horizontal boxes. Then the degree of the $\tau_k$ as a polynomial
of $t_1$ is $k\times(l-k+1)$. The sum of those degrees gives the height $|2\rho|$ where
$\rho$ represents the sum of all the fundamental weights $\omega_j$ \cite{flaschka:91}.

Then the geometry of the divisor ${\mathcal D}_0=\cup_{j=1}^l {\mathcal D}_j$, the union of the Painlev\'e
divisors ${\mathcal D}_j=\{\tau_j=0\}$ near the point $p_o$, can be expressed as the product of
$\tau$-functions,
\begin{equation}
\label{divF}
F(t_1,\ldots,t_l):=\prod_{j=1}^l\tau_j(t_1,\ldots,t_l)=F_d(t_1,\ldots,t_l)+F_{d+1}(t_1,\ldots,t_l)+\cdots\,,
\end{equation}
where each $F_k(t_1,\ldots,t_l)$ is a degree $k$ homogeneous polynomial of $(t_1,\ldots,t_l)$.
The algebraic variety $V:=\{F_d=0\}$ defines the {\it tangent cone} at the point $p_o$, and the degree
$d$ is the {\it multiplicity} of the singularity of $V$ at $p_o$. 
The degree $d$ is also related to the total number of blow-ups along the trajectory of the Toda flow
(Definition \ref{pq}).
In this article, we compute the number $d$ and express it as
 the number of ${\mathbb F}_q$-points of a finite Chevalley group (Proposition \ref{totalBnumber}); or equivalently as the dimension of any Borel subalgebra of ${\rm Lie}( \check K ({\mathbb C}))$
 (see \cite{carter}).
\begin{Example}
For $A_l$-Toda lattice, counting the {\it minimal} degree of the Schur polynomial $S_{(l-k+1,\ldots,l)}$
for each $\tau_k$,
one can easily find from (\ref{divF}) that $d=l(l+2)/4$ for $l$ even, and $d=(l+1)^2/4$ for $l$ odd.
The number $d$ agrees with  the total number of blow-ups in this Toda flow (\cite{kodama:98},
see also Section \ref{pq:section}).
The $d$ also gives the degree of the polynomial $\tilde p(q)=q^{-r}|K({\mathbb F}_q)|$ (see Propositions  \ref{direct} and \ref{totalBnumber}). Note in the case of $A_l$,
$K=SO(l+1)$.
\end{Example}


\subsubsection{Action of the Weyl group on the signs.}
Here we give an algebraic description of the blow-ups, so that one can compute
the number of blow-ups in the Toda flow:
The following action of the Weyl group $W$ on signs describes how the signs of the functions
$a_j$ for $j=1,\ldots, l$ change when $a_i$ blows up.

\begin{Definition} \label{act} (Proposition 3.16 in \cite{casian:02}) For any set of signs
$\epsilon =(\epsilon_1\ldots  \epsilon_l)\in \{\pm\}^l$, we define $W$-action on the signs
as follows:
A simple reflection $s_i:=s_{\alpha_i}\in W$ acts on the sign $\epsilon_j$ by
\begin{equation}
\nonumber
s_i~:~\epsilon_j \longmapsto \epsilon_j\epsilon_i^{C_{j,i}}.
\end{equation}
We also define the relation $\Rightarrow$ between the vertices of the polytope
$\Gamma_{\epsilon}$ as follows: if  $\epsilon_i=+$ then we  write 
$\epsilon \,\overset{s_i}{\Longrightarrow}\,\epsilon^\prime$ where  $\epsilon^\prime=s_i\epsilon$. We
also write $w\Rightarrow ws_i$ (see Figure \ref{A1:fig}).
\end{Definition}
Thus the sign change is defined on the group character $\chi_{\alpha_i}$ with
$\epsilon_i={\rm sgn}(\chi_{\alpha_i})$ (recall $s_i\cdot\alpha_j=\alpha_j-C_{j,i}\alpha_i)$. 
We here identify the sign $\epsilon_i$
as that of $a_i$, since the functions $a_k$ are given by (\ref{ab}) and the $\tau$-functions are
given by the fundamental weights $\omega_j$ in (\ref{tau}), which relate to the equation
$\chi_{\alpha_j}=\prod_{k=1}^l(\chi_{\omega_k})^{C_{j,k}}$ (recall $\alpha_j=\sum_{k=1}^l C_{j,k}\omega_k$).

Under the action of the Weyl group, not every simple reflection $s_i$ changes the
sign $\epsilon$. The following is an alternative way to measure the size of $w$
which only takes into account simple reflections that change the sign $\epsilon$,
that is, a trajectory of a Toda lattice having a blow-up point.
These numbers will later reappear in the context of the computation of certain
Frobenius eigenvalues.

Now the following definition gives the number of blow-ups in the Toda orbit from
the top vertex $e$ to the vertex labeled by $w\in W$:

\begin{Definition} \label{eta} 
For a reduced expression $w=s_{j_1}\cdots s_{j_r}$ and $\epsilon\in\{\pm\}^l$,
we define the number,
\[
\eta(w,\epsilon):=\left|\{~j_k~|~(s_{j_{k-1}}\cdots s_{j_1}\epsilon)_{j_k}=-\,,~k=1,\ldots,r~\}\right|\,,
\]
Namely, for the sequence of signs as the orbit given by $w$-action, 
\[ \epsilon~ \to ~s_{j_1}\epsilon~ \to ~s_{j_2} s_{j_1} \epsilon ~\to~ \cdots~\to~w^{-1}\epsilon\,.
\]
$\eta(w,\epsilon)$ is the number of $\to$ which are {\it not} of the form $\overset{s_i}{\Longrightarrow}$ as in Definition \ref{act}.
For example, if $\epsilon=(+\ldots+)$, then $\eta(w,\epsilon)=0,~\forall w\in W$.
The number $\eta(w_*,\epsilon)$ for the longest element $w_*$ gives the total
number of blow-ups along the Toda flow in $\Gamma_{\epsilon}$-polytope.
 Whenever $\epsilon=(-\ldots -)$ we 
will just denote $\eta(w,\epsilon)=\eta(w)$.
\end{Definition}

Note that each reduced expression of $w$ corresponds to a path following  Toda lattice
trajectories along 1-dimensional subsystems leading to $w$.  Each 1-dimensional
subsystem is equivalent to $A_1$-Toda lattice discussed in Introduction.
It can be shown that
$\eta(w,\epsilon)$ is independent of the reduced expression (Corollary \ref{indep}).  Hence the number of blow-up
points along trajectories in one-dimensional subsystems in the boundary of
the $\Gamma_{\epsilon}$-polytope is independent of the trajectory (parametrized
by a reduced expression).  
 
 We will find below a computation of  $\eta (w_*)$, the total number of blow-ups
 along the Toda flow, in terms of the polynomial computing the order of $K({\mathbb F}_q)$. Alternatively, $q^{-\eta(w_*)}$ is given in terms of the Hecke algebra operators of Lusztig and Vogan in \cite{lusztig83}  (Remark \ref {zz} combined with Proposition \ref{equalweights}). The expression   $q^{\eta(w_*)}$ is also the Frobenius eigenvalue occurring in the top (etale)  cohomology group of  the open dense $K(k_q)$ orbit  ${\mathcal O}_o(k_q)$  (Proposition \ref{frob}). 
 
 \begin{figure}[t]
\includegraphics[width=4.7in]{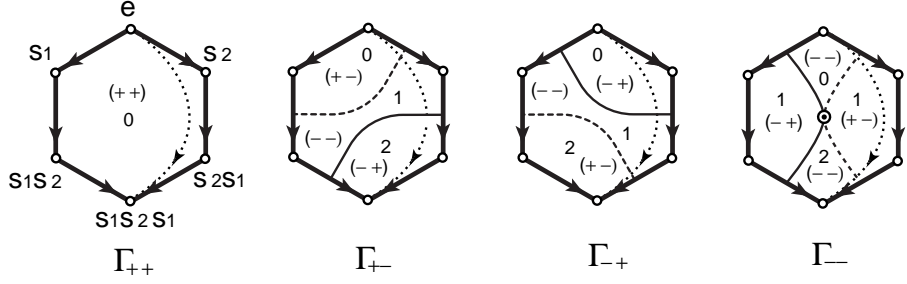}
\caption{The four hexagons $\Gamma_{\epsilon}$ associated to the $A_2$-Toda lattice. 
In each $\Gamma_{\epsilon}$, the Painlev\'e divisors ${\mathcal D}_j$
 are indicated by the solid curves for the blow-ups of $a_1$ (i.e. $\tau_1=0$)
and the dashed curves for the blow-up of  $a_2$ (i.e. $\tau_2=0$).
The boundaries of the hexagons describe the subsystems given by $a_i=0$ for $i=1,2$.
  Each hexagon is divided by the 
Painlev\'e divisors into connected components, e.g.
$\Gamma_{+-},\Gamma_{-+}$ have two, and $\Gamma_{--}$ has
four connected components. The Toda flow in $t_1$-variable
is shown as the dotted curve starting from
the vertex marked by the identity
element $e$, and ending to the vertex by the longest element $w_*=s_1s_2s_1$.
The numbers in the polytopes indicate $\eta(w,\epsilon)$, the number of blow-ups
along the Toda flow. The signs in each connected components in $\Gamma_{\epsilon}$
represent the signs $({\rm sgn}(a_1),{\rm sgn}(a_2))$, while $\epsilon$ in $\Gamma_{\epsilon}$
denotes the signs of $(a_1,a_2)$ in the connected component including the top vertex $e$.}
\label{hexagon1:fig}
\end{figure}

\begin{figure}[t]
\includegraphics[width=5.2in]{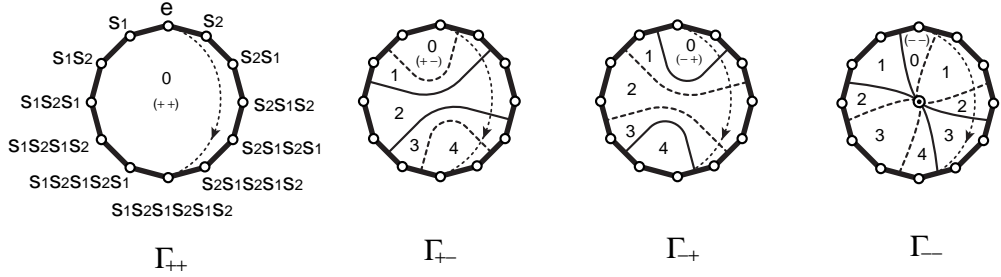}
\caption{The four 12-gons $\Gamma_{\epsilon}$ associated to the $G_2$-Toda lattice.
The Painlev\'e divisors corresponding to the blow ups of $a_1$ and $a_2$ are 
indicated by the solid and dashed curves, respectively. The Toda flow in $t_1$-variable
is shown as the dotted curve starting from
the vertex marked by $e$, and ending to the vertex by $w_*=s_1s_2s_1s_2s_1s_2$.
The numbers in $\Gamma_{\epsilon}$ represent $\eta(w,\epsilon)$.}
\label{G2:fig}
\end{figure}

\begin{Example}\label{A2G2}
The cases of $A_2$ and $G_2$ are illustrated in Figures \ref{hexagon1:fig} and \ref{G2:fig}.
In these figures the four hexagons and 12-gons are shown as the $\Gamma_{\epsilon}$-polytope with the signs $\epsilon=(\epsilon_1\epsilon_2)$.
These polytopes glue together to form a compact isospectral manifold $\tilde Z(\gamma)_{\mathbb R}
=\cup_{\epsilon\in\{\pm\}^2}\Gamma_{\epsilon}$. In \cite{casian:02b}, we have shown for example that $\tilde Z(\gamma)_{\mathbb R}$
gives a connected sum of two Klein bottles for $A_2$, and a connected sum of five Klein bottles for $G_2$.
Trajectories of the Toda lattice starts in the vertex associated to $e$ and move towards the vertex corresponding to the longest element $w_*$ in the Weyl group. 

In the case of $A_2$, the $W$-action on the signs $\epsilon=(\epsilon_1\epsilon_2)$ gives
  $s_1(--)=(-+)$,
$s_2(-+)=(-+)$ and $s_1(-+)=(--)$. From those we obtain $\eta (e)=0$, $\eta(s_1)=\eta(s_2)=\eta(s_1s_2)=\eta(s_2s_1)=1$ and $\eta(s_1s_2s_1)=2$. Those give the numbers of blow-ups
in the Toda flow (see Figure \ref{hexagon1:fig}).

In the case of $G_2$, we obtain $\eta(e)=0$, $\eta(s_1)=\eta(s_2)=\eta(s_1s_2)=\eta(s_2s_1)=1$,
$\eta(s_1s_2s_1)=\eta(s_2s_1s_2)=2$, $\eta(s_1s_2s_1s_2)=\eta(s_2s_1s_2s_1)=\eta(s_1s_2s_1s_2s_1)=\eta(s_2s_1s_2s_1s_2)=3$
and $\eta(w_*)=4$. The total number of blow-ups $\eta(w_*)$ is then 4 (see Figure \ref{G2:fig}).

\end{Example}

Before closing this Section, we mention a Poincar\'e duality relation satisfied by  the numbers $\eta(w)$:
\begin{Proposition}\label{poincare}
The number $\eta(w)$ satisfies the following  Poincar\'e duality relation,
\[
\eta(w_*)=\eta(w_*w)+\eta(w)\,.
\]
\end{Proposition}
\begin{Proof}
First note that $w_*(-\ldots -)=(-\ldots -)$. Also we note that $w_*w$ gives a path from $w_*$ to $w_*w$ which is dual to the path from $e$ to $w$ (recall $w$ acts on $\epsilon$ as $w^-1\epsilon$).
Since $(w_*w)^{-1}(-\ldots -)=w^{-1}(-\ldots -)$, we get
$\eta(w_*w)=\eta(w_*)-\eta(w)$.
\end{Proof}

\section{Blow-ups and graph of incidence numbers}

The strategy in the proof of our main results of constructing a partial dictionary 
between the objects in the Toda lattice and the flag manifold has the following order:

\begin{enumerate} 
 \item{} First define the polynomial $p(q)$ in terms of the blow-ups along the trajectory of the Toda flow (Definition \ref{pq}).
 \item{} Compute directly the multiplicity $d$ of the singularity of ${\mathcal D}_0=\cup_{j=1}^l{\mathcal D}_j$
 of the union of the Painlev\'e divisors at the point of intersection $p_o$ of the divisors, $\{p_o\}=\cap_{j=1}^l{\mathcal D}_j$. Then show $d$ is also given by the degree of the polynomial
 $\tilde p(q):=q^{-r}|K({\mathbb F}_q)|$, which is the dimension of any Borel subalgebra of ${\rm Lie}(K({\mathbb C}))$ (Proposition \ref{direct}).
 \item{} Define the graph ${\mathcal G}_{\epsilon}$ associated with the blow-ups in the polytope ${\Gamma}_{\epsilon}$ (Definition \ref{graph}). The signs $\epsilon_i:={\rm sgn}(a_i)$
 on $\Gamma_{\epsilon}$ are interpreted in terms of local systems arising in the computation of
 integral cohomology of the real flag manifold (Definitions \ref{thesigns}).
 \item{} Also define the graph ${\mathcal G}_{\mathcal L,\mathcal B}$ 
  in terms of a locally constant sheaf ${\mathcal L}$, which leads to $H^*({\mathcal B};{\mathcal L})$
  (Definition \ref{graphL}).
\item{} Show that the graphs ${\mathcal G}_{\epsilon}$ and
  ${\mathcal G}_{\mathcal L,\mathcal B}$ coincide (Theorem \ref{mainT} and Proof is given in Section \ref{mainTproof}). This completes the relation between the integral
  cohomology of the real flag manifold and the blow-up structure of the Toda lattice.
\item{} Count the number of blow-ups along trajectories of the Toda lattice, and relate it
to Frobenius eigenvalues in the flag manifold over positive characteristic (Proposition \ref{frob}).
\item{} Show $H^*({\mathcal B}; {\mathbb Q}) = H^*(K; {\mathbb Q})$ (Proposition \ref{flagvsgroup}).
\item{} Apply  Lefschetz fixed point theorem  to the Frobenius map to derive a relation
between the number of blow-up points along trajectories of the Toda lattice and  
the order of finite Chevalley group $K({\mathbb F}_q)$ (Theorem \ref{Kpq}).
\item{} Show that the degree of $p(q)$ gives the multiplicitiy $d$ of  the union of all the  Painlev\'e divisors at 
the point of intersection $p_o$, i.e. $d=\eta(w_*)$  (Proposition \ref{totalBnumber}).
\end{enumerate}

\subsection{The polynomial $p(q)$}\label{pq:section}
We now introduce polynomials in terms of the numbers $\eta(w,\epsilon)$
in Definition, \ref{eta} which play a fundamental role in this paper. We
refer to these polynomials as the {\it alternating sum of the blow-ups}.
For each $\Gamma_{\epsilon}$-polytope we then define $p_\epsilon(q)$:

\begin{Definition} We define a polynomial, the alternating sum of the blow-ups,
\begin{equation}
\label{pq}
 p_\epsilon(q)=(-1)^{l(w_*)}\sum_{w\in W} (-1)^{l(w)}q^{\eta(w,\epsilon)}\,,
\end{equation}
where $w_*$  is the
longest element of the Weyl group and $l(w)$ indicates the length of $w$. When $\epsilon=(-\ldots -)$, we 
simply denote this polynomial by $p(q)$. The degree of $p(q)$ gives
the total number of blow-ups along the Toda flow, $\eta(w_*)={\rm deg}(p(q))$.
\end{Definition} 

We then have:
\begin{Proposition}\label{pzero}
We have  $p_\epsilon (q)=0$ unless $\epsilon=(-\ldots -)$. 
\end{Proposition}
\begin{Proof}
If the $i$th component $\epsilon_i$ of $\epsilon$ is $+$, then we have
$s_i\epsilon=\epsilon$ and hence $\eta(e,\epsilon)=\eta(s_i,\epsilon)=0$.
Now for any $w\in W$ such that 
$l(s_iw)=l(w)+1$, we have  $\eta(w,\epsilon)=\eta(s_iw,\epsilon)$ (note that  $\epsilon = s_i\epsilon$ and applying $w^{-1}$ on both sides we also have 
$w^{-1}\epsilon = w^{-1}s_i\epsilon= (s_iw)^{-1}\epsilon$).  This implies 
$\eta(w,\epsilon)=\eta(s_iw,\epsilon)$. Then  for any $w$ of minimal length in its coset in  $\langle s_i \rangle\backslash W$,    the sum contains a pair $q^{\eta(w,\epsilon)}$ and
$q^{\eta(s_iw,\epsilon)}=q^{\eta(w,\epsilon)}$ of opposite signs, and hence they do not contribute the sum. By writing the sum over $W$ as a sum over the set of disjoint cosets of $\langle s_i \rangle\backslash W$  we obtain  $p_\epsilon (q)=0$.
\end{Proof}

Hence the only relevant
polynomial here is $p(q)$, the polynomial for the $\Gamma_{-}$-polytope. 
This proposition corresponds under our dictionary to the fact that the {\it rational}
cohomology of $K$ and ${\mathcal B}=K/T$ actually agree. 
 Recall that we are dealing only with the case
when the Lie algebra is split and the group $T$ is then a finite group.
The polynomial $p_\epsilon (q)$ corresponds, under our dictionary, to a Lefschetz number
for the Frobenius action over a field of positive characteristic 
for cohomology with local coefficients.  When $q=1$ these polynomials all vanish,
including $p(q)$. This reflects the fact that the Euler characteristic of $K$ is zero.

As we will explain below that the numbers $\eta(w,\epsilon)$ are deeply tied up with
the cohomology of the flag manifold. The polynomial $p(q)$ under our dictionary, contain
all the information regarding the cohomology ring of the compact
Lie group $K$. 

\begin{Example} From Figures \ref{hexagon1:fig}
and \ref{G2:fig}, we note that the numbers $\eta(w,\epsilon)$ are constant
on the connected components in a given $\Gamma_{\epsilon}$-polytope.  The polynomials $p(q)$
that are obtained from (\ref{pq}) by counting $\eta(w)$ are $p(q)=q^2-1$ for type $A_2$, and
$(q^2-1)^2$ for type $G_2$.  The $A_3$ example is a bit more complicated but
it gives the same polynomial $p(q)=(q^2-1)^2$ obtained in the $G_2$ case (see Figure \ref{A3a:fig} below).
The reason given earlier is that $K$ is in both cases essentially $SU(2)\times SU(2)$.
This will become clear when the connection with $K$ is made below.
\end{Example}

We are also interested in the degree of the polynomial $p(q)$, i.e. $\eta(w_*)$. 
One reason for this is the following
(see also below Eq. (\ref{divF})):

\begin{Conjecture}\label{tangentcone}
The degree $\eta(w_*)$ of the polynomial $p(q)$ is the multiplicity $d$ of the singularity of the
divisor given by ${\mathcal D}_0=\cup_{j=1}^l{\mathcal D}_j$ at the point $p_o$
where all the divisors ${\mathcal D}_j=\{\tau_j=0\}$ intersect. Namely the number $d$ is
given by the minimal degree of the product of the $\tau$-functions, i.e. (\ref{divF}).
Furthermore, the degree $\eta(w_*)$
also gives the number of {\it real} roots of the Schur polynomials associated to the nilpotent Toda lattices
(see \cite{casian:04} for the nilpotent Toda lattices).
\end{Conjecture}
The first part of the conjecture can be verified directly for the Toda lattice associated with any Lie algebra of the classical type or type $G_2$
by counting the minimal degrees of Schur polynomials,
which are the leading terms of the $\tau$-functions near the point $p_o$: First note that the $\tau$-functions for the nilpotent
Toda lattices are obtained from (\ref{tau}). Then find explicit forms of the highest weight vectors and
the companion matrix (the regular nilpotent element) for each algebra. The following is the result
based on the nilpotent Toda lattices:
\smallskip

 In the case of $A_l$-Toda lattice, the $\tau$-functions
are given by the Schur polynomials in (\ref{schur-tau}), i.e.
\[
\tau_{k}(t_1,\ldots,t_l)=(-1)^{\frac{k(k-1)}{2}}S_{(l-k+1,\ldots,l)}(t_1,\ldots,t_l)\quad k=1,\ldots,l\,.
\]
The minimal degrees of those Schur polynomials are given as follows:
\begin{itemize}
\item{} For $l$ even, the minimal degrees of $\tau_j$, $j=1,\ldots,l$, are given by
\[
1,~2,~\ldots~,~\frac{l}{2},~\frac{l}{2},~,\ldots~,~2,~1,
\]
(e.g. $\tau_1\sim t_l, ~\tau_2\sim t_{l-1}^2$ and $\tau_l\sim t_l$).
 The sum of those degrees then gives $d=\frac{l(l+2)}{4}$.
\item{} For $l$ odd, the minimal degrees are
\[
1,~2,~\ldots~,~\frac{l-1}{2},~\frac{l+1}{2},~\frac{l-1}{2},~\ldots~,~2,~1.
\]
 from which we have $d=\left(\frac{l+1}{2}\right)^2$.
\end{itemize}
\smallskip

In the case of $B_l$-Toda lattice, the $\tau$-functions are given by
\[
\tau_k(t_1,t_3,\ldots,t_{2l-1})={\rm Wr}(h_{2l},\ldots,h_{2l-k+1})\, \quad k=1,\ldots,l-1
\]
and 
\[\tau_l(t_1,t_2,\ldots,t_{2l-1})=\sqrt{|{\rm Wr}(h_{2l},\ldots,h_{l+1})|}\,,
\]
where $h_k$ are given in (\ref{schurp}) with $t_{2k}=0$ for all even parameters (see \cite{casian:04}).
Then we have:
\begin{itemize}
\item{} For $l$ even, the minimal degrees are given by
\[
2,~2,~4,~4,~\ldots~,~l-2,~l-2,~l,~\frac{l}{2}\,.
\]
(e.g. $\tau_1\sim t_1t_{2l-1},~\tau_2\sim t_{2l-1}^2$ and $\tau_l\sim (t_{l+1})^{l/2}$).
The degree $d$ is then given by $d=\frac{l(l+1)}{2}$.
\item{} For $l$ odd, the minimal degrees are
\[
2,~2,~4,~4,~\ldots~,~l-1,~l-1,~\frac{l+1}{2}\,,
\]
 from which we have $d=\frac{l(l+1)}{2}$.
\end{itemize}
\smallskip

In the case of $C_l$-Toda lattice, the $\tau$-functions are
\[
\tau_k(t_1,t_3,\ldots,t_{2l-1})={\rm Wr}(h_{2l-1},\ldots,h_{2l-k})\,\quad k=1,\ldots,l\,.
\]
Again one takes $t_{2k}=0$ for $h_n$. Then the minimal degrees are given by
\[
1,~2,~3,~\ldots~,~l-1,~l\,.
\]
This gives $d=\frac{l(l+1)}{2}$ (which is the same as $B_l$-case).
\smallskip

In the case of $D_l$-Toda lattice, the $\tau$-functions are given as follows:
\begin{itemize}
\item{} For $l$ even, they are given by, for $k=1,\ldots,l-2$,
\[
\tau_k(t_1,t_3,\ldots,t_{2l-3},s)={\rm Wr}(sh_{l-1}+2h_{2l-2},sh_{l-2}+2h_{2l-3},\ldots,sh_{l-k}+2h_{2l-1-k})\,.
\]
The $\tau_{l-1}$ and $\tau_l$ are given by
\[
[\tau_{l-1}\cdot\tau_l](t_1,t_3,\ldots,t_{2l-3},s)={\rm Wr}(sh_{l-1}+2h_{2l-2},\ldots,sh_{2}+2h_l)\,,
\]
and
\[
(\tau_l(t_1,t_3,\ldots,t_{2l-3},s))^2=\left|\begin{matrix}
sh_{l-1}+2h_{2l-2} & sh_{l-2}+2h_{2l-3}&\cdots &sh_{1}+2h_{l}&s+h_{l-1}\\
sh_{l-2}+2h_{2l-3} &sh_{l-3}+2h_{2l-4}&\cdots & s+2h_{l-1}  &  h_{l-2} \\
\vdots&\vdots  & \ddots  &  \vdots &\vdots  \\
sh_1+2h_l & s+2h_{l-1} &\cdots & 2h_2 & h_1\\
s+h_{l-1} & h_{l-2} &\cdots & h_1 & 0
\end{matrix}\right|
\]
Here the even parameters are all zero $t_{2k}=0$, and $s$ is a flow parameter
associated with the Chevalley invariant with degree $l$.
Counting the minimal degrees of the $\tau$-functions, we have
\[
2,~2,~4,~4,~\ldots~,~l-2,~l-2,~\frac{l}{2},~\frac{l}{2}\,,
\]
(e.g. $\tau_1\sim st_{l-1},~\tau_2\sim t_{2l-3}^2$ and $\tau_l\sim s^{l/2}$). Then we have $d=\frac{l^2}{2}$.
\item{} For $l$ odd, the $\tau$-functions are given by, for $k=1,\ldots,l-2$,
\[
\tau_k(t_1,t_3,\ldots,t_{2l-3},s)={\rm Wr}(s^2+2h_{2l-2},2h_{2l-3},\ldots,2h_{2l-1-k})\,.
\]
The last two $\tau$-functions are 
\[
[\tau_{l-1}\cdot\tau_l](t_1,t_3,\ldots,t_{2l-3},s)={\rm Wr}(s^2+2h_{2l-2},2h_{2l-3},\ldots,2h_{l})\,,
\]
and 
\[
(\tau_l(t_1,t_3,\ldots,t_{2l-3},s))^2=\left|\begin{matrix}
s^2+2h_{2l-2} & 2h_{2l-3}&\cdots &2h_{l}&s+h_{l-1}\\
2h_{2l-3} &2h_{2l-4}&\cdots & 2h_{l-1}  &  h_{l-2} \\
\vdots&\vdots  & \ddots  &  \vdots &\vdots  \\
2h_l & 2h_{l-1} &\cdots & 2h_2 & h_1\\
s+h_{l-1} & h_{l-2} &\cdots & h_1 & 1
\end{matrix}\right|
\]
Now the minimal degrees of the $\tau$-functions are
\[
2,~2,~\ldots~,~l-3,~l-3,~l-1,~\frac{l-1}{2},~\frac{l-1}{2}\,.
\]
Then we have $d=\frac{l^2-1}{2}$.
\end{itemize}
\smallskip

In the case of $G_2$, we have
\[
\tau_1(t_1,t_5)=h_6,\quad \tau_2(t_1,t_5)={\rm Wr}(h_6,h_5)\,.
\]
Here the parameters are only $t_1$ and $t_5$ and all others take zero.
The minimal degree is then 2 for each $\tau$-function, i.e.
$\tau_1\sim t_1t_5$ and $\tau_2\sim t_5^2$, and we have 
$d=4$.
\smallskip

We can summarize these computations in the following:

\begin{Proposition} \label{direct} Let  ${\mathfrak g}$ be a split semisimple Lie algebra not containing factors of type $E$ or $F$. The multiplicity $d$ of the singularity at $p_o$, given as the minimal degrees
of Schur polynomials  in $(\ref{schur-tau})$,  is the dimension of any Borel subalgebra of
${\rm Lie}( \check K  ({\mathbb C}))$. Moreover, $d$ is also the degree of the polynomial
$\tilde p(q)=q^{-r}|\check K({\mathbb F}_q)|$. 
\end{Proposition}
\begin {Proof} The statement concerning the dimension of any Borel subalgebra follows case-by-case from the explicit computations that are listed above.  This relies on an approximation near $p_o$  of the tau functions in the semisimple case by the tau functions in the nilpotent case (given by Schur polynomials).  One  also needs the list of groups $K$ that appear in each case. This is provided, for example,  in the satement of  Theorem \ref {Kpq}.  The second statement uses a formula for the sum of the degrees of basic invariant polynomials. This sum gives the degree of 
the polynomial $\tilde p(q)=q^{-r}|\check K({\mathbb F}_q)|$.  This formula is Theorem 9.3.4 of \cite{carter}. 
\end{Proof}
In Proposition \ref{totalBnumber}, those numbers $d$ are shown to be also the degrees $\eta(w_*)$ of 
$p(q)$.  Note that $d$ is the multiplicity of a singularity and  $\eta(w_*)$ is defined differently, as the maximal number of blow-ups encountered along the Toda flow, counted along one dimensional subsystems. These polynomials $p(q)$, the alternating sum of the blow-ups, are then shown to agree with  the polynomials  $\tilde p(q)$. Also notice that the minimal degree
for each $\tau$-functions is quite similar to each degree $d_i$ of the basic $W$-invariant
polynomial of the Chevalley group $K$ (see \cite{carter} or Theorem \ref{Kpq} below). We did not compute  the exact relation between those
degrees, but we expect that each degree $d_i$ is related to the number of real intersection points
on the tangent cone $V=\{F_d=0\}$ defined in (\ref{divF}) with a linear line
corresponding to the $t_1$-flow of the Toda lattice. This may be stated as
\[
{\rm deg}(p(q))=\underset{c\in{\mathbb R}^{l-1}}{\rm Max}\left|\{{\mathcal D}_0\cap L_c~|~
{\rm transversal ~intersection}\,\}\right|
\]
where ${\mathcal D}_0=\cup_{j=1}^l\{\tau_j=0\}$ and $L_c:=\{(t_1,c_2,\ldots,c_l)\,|\,c=(c_2,\ldots,c_l)\in{\mathbb R}^{l-1}\}$. It is interesting to note that the
degrees $d_i$ and the minimal degrees of $\tau$-functions are the same for the cases having
the same ranks, $l={\rm rank}(\mathfrak g)={\rm rank}(K)$, i.e.
the cases of $B,\,C$ and $D_l$ with $l$ even  (see Theorem \ref{Kpq}).

\subsection{The graphs associated to the blow-ups}

The following graph ${\mathcal G}_\epsilon$ was originally  motivated by the problem of computing
the number of connected components in the $\Gamma_{\epsilon}$-polytope. This problem
is analogous to the problem of computing the intersection of two opposite 
top dimensional Bruhat cells in the case of a real flag manifold (e.g. see \cite{shapiro:97,rietsch:97,zelevinsky:00}). 
Here we were counting the number of connected components appearing
in the intersection $(N^*B/B)\cap(G^{C_{\gamma}}w_*B/B)$ where
$G^{C_{\gamma}}$ is defined in (\ref{gcorbit}) (note in particular,
when $\gamma=0$ (nilpotent case), $G^{C_0}\subset N$). Also note that the Painlev\'e
divisor ${\mathcal D}_J:=\cap_{j\in J}{\mathcal D}_j$ for $J\subset\{1,\ldots,l\}$ corresponds to
the $N^*$-orbit, $N^*w_JB/B$, where $w_J$ is the longest element of the subgroup
$W_J:=\langle s_i\,|\,i\in J\rangle$ (see \cite{flaschka:91,casian:02b}).
We then observed that in all examples
this was the graph of incidence numbers for a real flag manifold, observation
which then started the present study.

\begin{Definition}
\label{graph}
For a fixed $\epsilon=(\epsilon_1\ldots \epsilon_l)$, we associate a graph
${\mathcal G}_\epsilon$ to the blow-ups of the
Toda lattice. The graph consists of vertices labeled by the elements
 of the Weyl group $W$, i.e. the vertices of the $\Gamma_{\epsilon}$-polytope, and 
 oriented edges $\Rightarrow$. The edges are defined as follows:
For any $w_1 , w_2 \in W$, there exist an edge between $w_1$ and $w_2$,
\[
w_1\Rightarrow w_2\quad {\rm iff}\quad
\left\{\begin{array}{llll}
{\rm a)}~ w_1\le w_2\,\,({\rm Bruhat~order})\\
{\rm b)}~l(w_2)=l(w_1)+1\,, \\
{\rm c)}~\eta(w_1,\epsilon)=\eta(w_2,\epsilon)\,,\\
{\rm d)}~w_1^{-1}\epsilon=w_2^{-1}\epsilon\,.
\end{array}\right.
\] 
When $\epsilon=(-\ldots -)$, we simply denote ${\mathcal G}={\mathcal G}_\epsilon$.
\end{Definition}

Note that  if there exists a path from $w_1$-vertex to
$w_2$-vertex without crossing a Painlev\'e divisor (i.e. no blow-up),
then we expect to have the edge, $w_1\Rightarrow w_2$. Namely, in this case,
$w_1\Rightarrow w_2$  means that the vertices associated to $w_1$
and $w_2$ belong to the {\it same} connected component of the polytope when blow-ups
are removed.  In particular this implies that the number of connected components, say $N_c$,  in the negative polytope $\Gamma_{-}$ is bounded below by the number of connected component, say $N_g$,  of the graph $\mathcal G$;  i.e.
we expect $N_g\le N_c$. For  type $A_l$ we have directly  verified  the equality
$N_c=N_g$ for    $l=1,2, 3, 4$  ($N_g=N_c$  in $A_3$ can be obtained from  Figures 4 and 5 ),  we expect the equality to be true for all $l$, and  in the negative polytope $\Gamma_{-}$.  However, for the case of $B_3$, we have $N_c=18$ and $N_g=17$. Finding the precise numbers of $N_c$ and $N_g$ is a very interesting problem. We can classify   the connected components in the negative polytope by dividing them into  families  according to  the signs of $a_i$ $i=1,\ldots , l$
(see Figures \ref{hexagon1:fig} and \ref{G2:fig}).  If we consider only {\it negative} components, that is,  components where  $a_i <0$ for all $i$, then computations suggest that there are exactly $2^g$ negative components with $g:={{\rm rank}( \check K)}$. This is the sum of all the Betti numbers of $\check K$. This formula works also in $B_3$ and \lq\lq explains\rq\rq the discrepancy between $N_c$ and $N_g$ in this case.

Summarizing: in  many cases the graph ${\mathcal G}_\epsilon$  accomplishes the job of joining together in its connected components exactly  those vertices $w$ of the polytope $\Gamma_\epsilon$
belonging to the same connected components.  This will be the case in the example considered below. 

\begin{Example} \label{exam}
In the case of $A_2$, we have $s_1(--)=(-+),\, s_2(-+)=(-+)$ which implies
$(--) \to (-+)\, {\overset{s_2}{\Longrightarrow} }\,(-+)$
and  $\eta(s_1)=1$, $\eta(s_1s_2)=1$. Therefore the graph ${\mathcal G}$
which encodes blow-up information in  $\Gamma_{--}$ of  Figure 2  is ($s_i $ is replaced with $i$):


$$\begin{matrix} 
{}                     &{e}& {}               & {} & {}& : & q^0 \\
[1]                     &{} & [2]                  & {} &{}& : & q^1 \\
 \Downarrow&{}&\Downarrow & {} & {} &{}& \\
  [12]               &{}& [21]                  & {} &{}& : & q^1\\
   {}&{[121]}&{ }                      &{}&{}& : & q^2\\
    \end{matrix} $$
    where the elements of the Weyl group are denoted by $s_{i}\cdots s_j=[i\cdots j]$. Here we have also  listed the monomials  $q^{\eta (w)}$ (in the variable $q$) associated to 
 representatives of the integral cohomology ($w\to \eta (w) \to q^{\eta (w)}$).
As already noted, the vertices of the hexagon $\Gamma_{--}$ belonging
to a connected component in Figure \ref {hexagon1:fig} form a connected
component of this graph (see also Figure \ref{A1:fig}). This graph classifying 
connected components in the hexagon minus the blow-ups agrees
with the graph in p.465 of  \cite{casian99} which is defined very differently  in terms of
incidence numbers. The graph of incidence numbers 
 gives rise to a chain complex by replacing the edges $\Rightarrow$ with
multiplication by $\pm 2$,
$${\mathbb Z}\langle e\rangle~ {\overset{\delta_0}{\longrightarrow} }~
{\mathbb Z}\langle s_1\rangle \oplus {\mathbb Z}\langle s_2\rangle ~{\overset{\delta_1}{\longrightarrow}}~ {\mathbb Z}\langle s_1s_2\rangle\oplus {\mathbb Z}\langle s_2s_1\rangle~ {\overset{\delta_2}{\longrightarrow}} ~{\mathbb Z}\langle s_1s_2s_1\rangle\,.  $$
Here $\langle w\rangle$ is the Bruhat cell associated to the element $w\in W$.
The only non-zero map is $\delta_1$ given by a diagonal matrix  with $2$'s in the diagonal corresponding to the $\Rightarrow$. The integral cohomology that results is
$$\left\{ \begin{matrix}
H^0(G/B;\mathbb Z )&= &\mathbb Z &  {} &: &{q^0}\\  
H^1(G/B;\mathbb Z )&=& 0                & {} &:& {q^1}\\
H^2(G/B;\mathbb Z )&=&\mathbb Z/2\mathbb Z \oplus\mathbb Z/2\mathbb Z & {} &:& {q^1}\\
H^3(G/B;\mathbb Z)&=&\mathbb Z & {} &:&  {q^2}\\
\end{matrix}\right. $$
Over the rationals this just gives the cohomology of $K=SO(3)$ which can be written as
the exterior ring with one generator $x_1$ of degree 2, i.e. $H^*(SO(3);{\mathbb Q})= 
\Lambda  (x_1q^2)$. The alternating sum of the 
$q^{\eta (w)}$ produces the polynomial $p(q)=q^2-1$. Now $p(q)$  multiplied by $q^{r}$ with
$r:={\rm dim}(K)-{\rm deg}(p(q))=1$ gives the number of points
of $SO(3)$ over a field with $q$ elements, i.e. $q^rp(q)=|SO(3;{\mathbb F}_q)|$  ($q$ is a power of an odd prime $p$, and
$r$ is the number of positive roots associated with the group $K$, \cite{carter}).
The explanation for this is   that the $q^{\eta (w) }$ listed
are  Frobenius eigenvalues in etale cohomology
of the appropriate varieties reduced to a field of positive characteristic (see p.465 of \cite{casian99}).
When this is taken into account we see that, in fact, more than
the cohomology of $K/T$ or $K$, we are obtaining  etale $\bar {\mathbb Q}_m$
cohomology over a field of  positive characteristic, including
Frobenius eigenvalues. All are derived from the structure of the Toda
lattice and its blow-up points.  This is analogous to the situation described  in Section \ref{coh} for the case of $A_3$.

If we start with $\epsilon=(-+)$ then we obtain the edges
$e\Rightarrow s_2$, $s_1\Rightarrow s_2s_1$, $s_1s_2\Rightarrow s_1s_2s_1$.
This time the vertices of the hexagon $\Gamma_{-+}$ belonging
to a connected component in Figure \ref {hexagon1:fig} form a connected
commponent of this graph. This graph now corresponds to the graph of incidence numbers computing
cohomology with local coefficients. The local system can be described 
by the signs $(-+)$. The $-$ sign indicates that along a circle in $G/B$
that corresponds to $s_1$, the local system is constant, and the second sign $+$ indicates
that along a circle corresponding to $s_2$, it is non-trivial. This 
can be made precise below and is an important step towards the dictionary.
The cohomology one obtains from this graph associated to $\Gamma_{-+}$  is the following:
$$\left\{\begin{matrix}
H^0(G/B;{ \mathcal L})&=&0  \\  
H^1(G/B; { \mathcal L})&=&\mathbb Z/2\mathbb Z   \\
H^2(G/B;{ \mathcal L} )&=&\mathbb Z/2\mathbb Z   \\
H^3(G/B;{ \mathcal L})&=&\mathbb Z/2\mathbb Z \\\end{matrix}\right.$$
\end{Example}
Here we did not list  the powers $q^{\eta(w,-+)}$ but clearly  $p_{-+}(q)=0$ (Proposition \ref{pzero}),
which implies that the rational cohomology with twisted
coefficients ${\mathcal L}$ is zero. Similarly we have the same results for $\Gamma_{+-}$ and $\Gamma_{++}$.

 \begin{figure}[t!]
\includegraphics[width=4.5in]{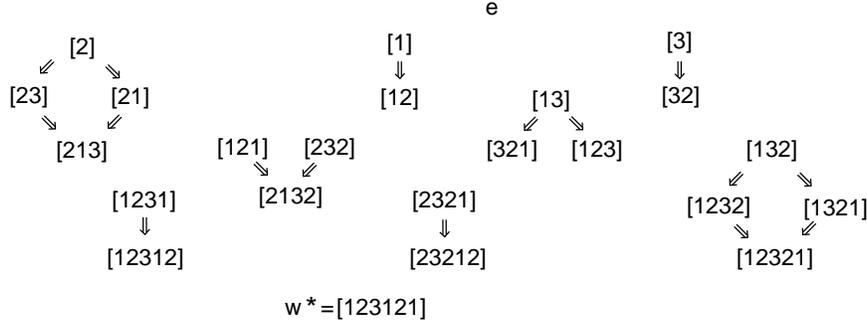}
\caption{The graph $\mathcal G$ of the real flag manifold for type $A_3$.
The Bruhat cells $NwB/B$ are denoted by $[ij\ldots k]$ for $w=s_is_j\ldots s_k$.
The incidence numbers associated with the edges $\Rightarrow$ are $\pm 2$ (see also
Example (8.1) in \cite{casian99}).
There are 10 connected components in this graph corresponding to the 10 connected components in $\Gamma_{-}$ after blow-up points are removed. }
\label{A3inc:fig}
\end{figure}

\begin{figure}[t!]
\includegraphics[width=4in]{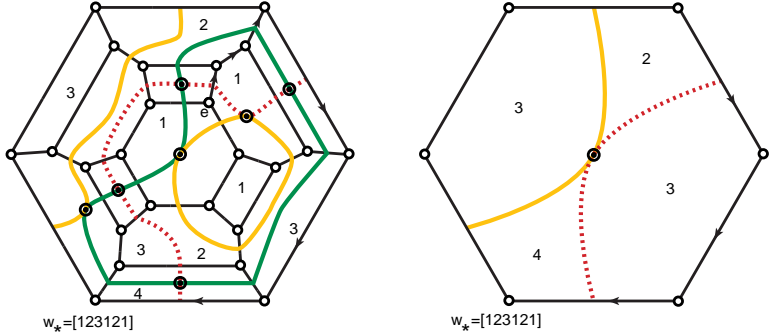}
\caption{The $\Gamma_-$-polytope for type $A_3$ and the Painlev\'e divisors (the right figure
is the back view of the left one). The
Painlev\'e divisors are shown by the dotted curve for ${\mathcal D}_1$, by the light color one for ${\mathcal D}_2$, and by the dark one for ${\mathcal D}_3$. The double circles indicate the divisor
${\mathcal D}_{ij}={\mathcal D}_i\cap{\mathcal D}_j$, which are all connected at the center of
the polytope $p_o$. The numbers $\eta(w)$ are obtained by using any path from $e$ to $w$  along edges of the polytope,  following the direction of the Toda flow, i.e. the path $w\to ws_i$ with
$l(ws_i)=l(w)+1$,  and counting intersections with the Painlev\'e divisors. The incidence graph in Figure \ref{A3inc:fig} can then be obtained from this Figure.}
\label{A3a:fig}
\end{figure}

 In the case of a Lie algebra of type $A_3$ we obtain the graph $\mathcal G$ in Figure \ref{A3inc:fig}.
This graph corresponds to the  polytope in Figure \ref{A3a:fig}  as separated into connected components by the divisors shown.  To determine the number $\eta(w)$ for any given $w$ it is enough to go from $e$ to $w$ along any path along the boundary   corresponding to a reduced expression $w=s_{n_1}\cdots s_{n_r}$ and count the number of intersections with the divisors. In Figure \ref{A3a:fig}, we show the path
following the expression $w_*=[123121]$, i.e. $e\to s_1\to s_1s_2\to s_1s_2s_3\to \cdots \to w_*$, with the arrows on the edges of the polytope.

We then state the following theorem showing the equivalence between the connected components
in the polytopes $\Gamma_{\epsilon}$ and the graphs ${\mathcal G}_{\epsilon}$ of the
incidence numbers defined in \cite{casian99}. A proof of the theorem will be given in Section \ref{mainTproof}.
\begin{Theorem}\label{mainT} The graph ${\mathcal G}$ is the graph of incidence numbers
for the integral cohomology of  the real flag manifold $\check {\mathcal B}$  in terms of 
the Bruhat cells.  In general, each $\check K$-equivariant local system  ${\mathcal L}$ on $\check {\mathcal B}$ corresponds to an $\epsilon$ and  ${\mathcal G}_\epsilon$ is the graph of incidence numbers for cohomology of $\check {\mathcal B}$ with twisted coefficients in ${\mathcal L}$.
\end{Theorem}
 
 \begin{Remark} In general there may be some signs $\epsilon$ such that no $\check K$-equivariant local system ${\mathcal L}$ corresponds to it.  For instance  if $G=SL(4;{\mathbb R})$ this is the case.  By considering the (possibly disconnected) group $\tilde G={\rm Ad}(SL(4;{\mathbb R})^{\pm})$ (see Remark \ref{group1}) this difficulty disappears;  this is mostly due to the fact that $H_{\mathbb R}$, the Cartan subgroup of  $\tilde G$,  has $2^l$ connected components which are not related by central elements. 
  \end{Remark}

\section{Real flag manifold $K/T$ and topology of $K$}

For each real split simple Lie algebra ${\mathfrak g}$  we have defined a connected Lie group $G$ and a maximal compact Lie subgroup $K$ (fixed point sets of a Cartan involution $\theta$). Moreover  in  the appropriate context of algebraic groups, all these objects can also be considered over a field $k$, an algebraic closure of a finite field ${\mathbb F}_q$ with $q$ elements.  We also recall  that the real flag manifold is first replaced with a complex manifold ${\mathcal O}_o$, a $K({\mathbb C})$-orbit in  the complex flag manifold ${\mathcal B}_{\mathbb C}$.

 
We now consider  the $K$-equivariant local systems on $G/B$ and ${\mathcal O}_o$  which arise from
the volume forms on the Bruhat cells. These  determine the incidence
numbers in a chain complex that computes integral cohomology of $G/B$
in \cite{casian99}.   The tangent space of the Bruhat cell associated to an element $w\in W$
is given by $w{\mathfrak n}^*\cap {\mathfrak n}$. The cotangent space
then corresponds to  $w{\mathfrak n}\cap  {\mathfrak n}^*$.  Now
 $\Lambda^{l(w)}w{\mathfrak n}\cap  {\mathfrak n}^*$ represents a volume form
on the Bruhat cell.  From here we can also obtain $l(w)$-forms on $G/B$
which are part of a de Rham chain complex.  
The action of $T$ on each  $\Lambda^{l(w)}w{\mathfrak n}\cap  {\mathfrak n}^*$
determines a $K$-equivariant local system. These local systems arising
from those  forms associated to Bruhat cells determine the  incidence
numbers with respect to (dual) of Bruhat cells.  Moreover the 
local systems determined by the  $\Lambda^{l(w)}w{\mathfrak n}\cap  {\mathfrak n}^*$
correspond to the signs $\epsilon_i={\rm sgn}(a_i)$ on the connected components of the
$\Gamma_{-}$-polytope. 

\begin{Example} Case of $SL(2;k)$\label{list1}. There are two $K(k)$-equivariant local systems on $k^*$ corresponding to the characters 
of a two element group $T=\left\{ \left(
\begin{matrix}
\epsilon_1 & 0   \\
0 & \epsilon_1 \\
\end{matrix}\right)\Big|  \epsilon_1=\pm 1\right\}$. The trivial character gives a constant sheaf ${\mathcal C}$  and
the non-trivial character gives a non-trivial local system ${\mathcal L}$.  Note that in this case $T$ acts trivially on $\Lambda^{l(w)}w{\mathfrak n}\cap {\mathfrak n}^*$.
For example if $w=s_1$, $\Lambda^{1}s_1{\mathfrak n}\cap {\mathfrak n}^*= {\mathbb R} X_{e_1-e_2}$,
 ($\alpha_1=e_1-e_2$)
and $T$ acts trivially on $X_{e_1-e_2}=\left( \begin{matrix}
0 & 1   \\
0 & 0 \\
\end{matrix}\right)$. 

The two local systems described correspond to the connected components 
in $\Gamma_{-}$ for the ${\mathfrak{sl}}(2; {\mathbb R})$ case of the example in Introduction. The non-trivial local system
corresponds to the single component $\Gamma_+$.  Tensoring with $\Lambda^{l(w)}w{\mathfrak n}\cap {\mathfrak n}^*$
produces the signs in $\Gamma_+$, we obtain for $w=s_1$, $\epsilon_1={\rm sgn}(a_1)=+$ corresponding to no blow-ups and to a non-trivial local systems in $G/B$. 
\end{Example} 
 
 We now make precise the correspondence between the various polytopes $\Gamma_\epsilon$
 and the computation of cohomology with local coefficients of $G/B$.

\subsection{Volume forms on the Bruhat cells and local systems}

In \cite{casian99} a chain complex of differential forms analogous to the de Rham chain
complex is used. The differential forms are sections of the line bundle
on the flag manifold $K/T$ induced by the action of $T$ on 
$\Lambda^{l(w)}w{\mathfrak n}\cap  {\mathfrak n}^*$.  The approach
in \cite{casian99} is to relate the boundaries in the chain complex to the representation
theory of the non-compact Lie group $G$ acting on these sections. Recall that a $K({\mathbb C})$-equivariant local system on ${\mathcal O}_o$
is given by a character of $T$.  If $H_+$ is a split Cartan subgroup in $G$  with connected component 
of the identity, then $H/H_{+}$ is isomorphic to $T$. Therefore,
since the Weyl group $W$ acts on both $H$ and $H_{+}$, it also acts on $T$.
From here it follows that $W$ acts on the characters of $T$. Given $\chi:T\to \{{\pm 1}\}$
character of $T$, we denote $w(\chi)$ the character obtained by applying $w$
to $\chi$. We will define a function $\epsilon$ associating to each $K$-equivariant local system ${\mathcal L}$ on $\check {\mathcal B}$ a list of $l$ signs $\epsilon$.    The action of $W$ on the characters of $T$ is just the action of $W$ on  signs which was introduced earlier (Definition \ref {eta}).

 \begin{Definition} Given the $K({\mathbb C})$-equivariant local system ${\mathcal L}$
corresponding to a character of $T$ given by $\chi({\mathcal L})$ we define
 ${\mathcal L}_w$ to be the $K({\mathbb C})$-equivariant local system
on ${\mathcal O}_o$ determined by the action of $T$ on
 $\chi({\mathcal L})\otimes \Lambda^{l(w)}w{\mathfrak n}\cap {\mathfrak n}^*$. 
\end{Definition}

We now consider a parabolic subgroup $P_i$ with complexification $P_i({\mathbb  C})$ associated to a simple root $\alpha_i$
containing the Borel subgroup $B$  of $G$.   There are projections:
$\pi_i: G({\mathbb C})/B({\mathbb C})\to G({\mathbb C})/P_i({\mathbb C})$ and
$\pi_i: G/B\to G/P_i $. Then restricting the local system ${\mathcal L}$
 to the fibers of $\pi_i$ can be described as follows:
Recall that the fibers of $\pi_i$ are ${\mathbb P}^1$, that is, the flag manifold
obtained from $SL(2;{\mathbb C})$. Each simple root $\alpha_i$ gives rise to a Lie group map,
$\Phi_{\alpha_i}: SL(2; {\mathbb R}) \to G$.  We consider $\Phi_{\alpha_i} ( \left(
\begin{matrix}
-1 & 0   \\
0 & -1 \\
\end{matrix}\right)) \in T$. Similarly we have maps $\Phi_{w\alpha_i}$ for $w\in W$.

\begin{Proposition} The local system  ${\mathcal L}_w$ is trivial along
the fiber of $\pi_i$ containing $x_w$ if and only if 
the character $\chi_{w}({\mathcal L})$  of $T$ on 
$w^{-1}\chi({\mathcal L})\otimes\Lambda^{l(w)}{\mathfrak n}\cap {w^{-1}\mathfrak n}^*$
is trivial on  $z_i=\Phi_{\alpha_i} ( \left(
\begin{matrix}
-1 & 0   \\
0 & -1 \\
\end{matrix}\right)).$

\end{Proposition}
\begin{Proof} 
 The local system  ${\mathcal L}_w$ is trivial along
the fiber of $\pi_i$containing $x_w$ if and only if the character  of $T$ on
$\chi({\mathcal L})\otimes\Lambda^{l(w)}w{\mathfrak n}\cap {\mathfrak n}^*$ is trivial on 
 $\Phi_{w\alpha_i} ( \left(
\begin{matrix}
-1 & 0   \\
0 & -1 \\
\end{matrix}\right)).$
Translating by  $w^{-1}$ everything  we obtain our statement.
\end{Proof}

We now recall that $T=K\cap B= K({\mathbb C})  \cap B({\mathbb C})$ is a finite group, and
the characters $\chi : T\to \{\pm\}$  parametrize $K({\mathbb C})$-equivariant local systems.   We will now find a more convenient  of expressing these  local systems in terms of $l$ signs which keep track of their structure along certain directions given by the simple roots. 

\begin{Definition} \label{thesigns} We let $\epsilon(w, {\mathcal L})= (-\chi_{w}({\mathcal L})(z_1),\ldots, -\chi_{w}({\mathcal L})(z_l))$.
If ${\mathcal L}$ is trivial we just write  $\epsilon(w, {\mathcal L})=\epsilon(w)$.
\end{Definition}

Note here that a constant sheaf ${\mathcal L }$ corresponds to a sign $\epsilon(w, {\mathcal L})=(-\ldots -)$.  With this parametrization we will have a correspondence in which a polytope $\Gamma_\epsilon$
corresponds to the choice of a local system ${\mathcal L}$  whose sign
 $\epsilon(e, {\mathcal L})=\epsilon$. The signs of the $a_i$ appearing inside the
polytope will agree with the various  $\epsilon(w, {\mathcal L})$ derived from the volume forms on the Bruhat cells for a fixed ${\mathcal L}$ (e.g. the trivial sheaf ${\mathcal L}$) in the dual  flag manifold 
 $\check {\mathcal B}$.  This is shown in Proposition  \ref{propsigns}. There we show that for fixed 
 ${\mathcal L}$    the set of signs $\{ \epsilon (w, {\mathcal L})\,|\, w\in W \}$ is a $W$-orbit. Since  the signs 
 $\{ ({\rm sgn}(a_1) \ldots {\rm sgn}(a_l) )\} $ corresponding to the polytope
  $\Gamma_{\epsilon (w, {\mathcal L})}$  are also given by the same $W$-orbit,  a correspondence is established between local systems associated to a fixed set of local coefficients ${\mathcal L}$,   and the signs appearing from the Toda lattice associated to a fixed polytope  
  $\Gamma_{\epsilon (w, {\mathcal L})}$. This is part of the \lq\lq partial dictionary\rq\rq.

\subsection{The graphs of incidence numbers of the real flag manifold}

We first review the computation of integral cohomology of  $G/B$ with certain $K$-equivariant local coefficients: 

Recall that there is  filtration by Bruhat cells with ${\mathcal B}_j:=\cup_{l(w)\le j}{NwB/B}$,
\[ 
\emptyset\, \subset \, {\mathcal B}_0\, \subset  \, {\mathcal B}_{1} \, \subset \, \cdots\, \subset \, {\mathcal B}_{l(w_*)}= {G/B}
\]
 Then there are coboundary maps 
$\delta:H^s({\mathcal B}_s, {\mathcal B}_{s-1}; {\mathbb Z}) \to 
H^{s+1}({\mathcal B}_{s+1}, {\mathcal B}_{s};{\mathbb Z})$ which
give rise to a chain complex computing the cohomology of $G/B$
(as in \cite{munkres} Theorem 39.4). We have $H^s({\mathcal B}_s, {\mathcal B}_{s-1};{\mathbb Z})={\rm Hom}_{\mathbb Z}( H_s({\mathcal B}_s, {\mathcal B}_{s-1};{\mathbb Z}),{\mathbb Z})$.
Also  $H_s({\mathcal B}_s, {\mathcal B}_{s-1};{\mathbb Z})=\oplus_{l(w)=s}H_s(\bar {\mathcal B}_w, \bar {\mathcal B}_w \setminus {\mathcal B}_w)$.
Denote $[w]$ the homology class that corresponds to ${\mathcal B}_w$. 
Generators of $H^{s+1}({\mathcal B}_{s+1}, {\mathcal B}_{s};{\mathbb Z})$  can 
 then be given by the ${\mathbb Z}$-module maps $f_w$ with $w\in W$ such that $l(w)=s$ 
determined  by $f_w([w'])=0$ if $w\not=w'$ and $f_w([w])=1$.
The coboundary can be written  in the form,
 $\delta(f_w)=\sum_{l(w')=s+1} [w;w']\,w'$. The incidence numbers are known
 to be $0$ or $\pm 2$ (see \cite{casian99}).

 \begin{Definition} \label{graphL}
 We define a graph ${\mathcal G}_{\mathcal B}$ whose
vertices consists of the elements in $W$ and oriented edges
$w\Rightarrow w'$ if and only if $[w;w']\not=0$, i.e. ${\mathcal G}_{\mathcal B}$
is the graph of incidence numbers associated with $H^*({\mathcal B};{\mathbb Z})$.
If  local coefficients given by a locally constant sheaf ${\mathcal L}$
 are used, then we denote the graph of incidence numbers leading to 
 $H^*({\mathcal B}; {\mathcal L})$ by ${\mathcal G}_{\mathcal L,B}$.
\end{Definition}

We now can rewrite Definition  7.4 b) of \cite{casian99} for the real split cases
in terms of these signs.  This will lead to the definition
of an oriented graph whose vertices are the elements of $W$.

\begin{Definition}\label{graphdef} We fix ${\mathcal L}$ and recall that
the sign associated to ${\mathcal L}_w$ is  
$\epsilon(w, {\mathcal L})=(\epsilon_1\ldots \epsilon_l)$.
We then define $w\, {\overset{\mathcal L} {\Longrightarrow } } \,ws_i$
 if and only if $\epsilon_i=+$, that is,
 the local system ${\mathcal L}_w$ is nontrivial along the fiber of $\pi$
 containining $x_w$. In addition, if $x$ increases the length of both $w$ and $ws_i$ and $w\,\overset{\mathcal L}{\Longrightarrow}\,ws_i$, then we also have an edge $wx\,\overset{\mathcal L}{\Longrightarrow}\,ws_ix$.
 \end{Definition}
 
 We now define the number $z(w,{\mathcal L})$, which is given in a way that is parallel to Definition \ref{eta} for the $\eta(w,\epsilon)$ computing the number of blow-ups along the Toda flow.
 However the $\epsilon$ keeping track of a polytope $\Gamma_{\epsilon}$ 
  is replaced with  $\epsilon(e,{\mathcal L})$ keeping track of a local system,  and each
 $\Rightarrow$ is replaced with ${\overset{\mathcal L}{\Longrightarrow }} $ i.e. a non-zero incidence number. Hence in $z(w,\mathcal L)$, we take into account not simple reflections that represent the crossing of a blow-up, but rather simple reflections that give a non-trivial coboundary.
 
 \begin{Definition} \label{zw}
 For any reduced expression $w=s_{j_1}\cdots s_{j_{l(w)}}$,
 we have $ s_{j_1}\cdots  s_{j_r} \to   s_{j_1}\cdots  s_{j_r} s_{j_{r+1}}$ for $r=1,\cdots,l(w)-1$.
 We let $z(w, {\mathcal L})$ denote the number of all $\to$ obtained as $r$ varies over $r=1,\ldots,l(w)-1$ which are {\it not}
 of the form ${\overset{\mathcal L}{\Longrightarrow}}$.  
 \end{Definition}
 \begin{Remark} \label{zz}
 The   numbers $z(w,\mathcal L)$ have the following  algebraic description (\cite {casian86,casian99}) . Consider the set 
${\mathcal D}$ of all $K({\mathbb C} )$-equivariant  local systems on all  the $K({\mathbb C} )$-orbits in 
 $G({\mathbb C})/B({\mathbb C} )$. Then as in  \cite{lusztig83}, the set  ${\mathcal M}$  of formal linear 
 combinations of elements in ${\mathcal D}$  with coefficients in ${\mathbb Z}[q, q^{-1}]$, ${\mathcal 
 M}={\mathbb Z}[q, q^{-1}]\otimes_{\mathbb Z} {\mathbb Z}[{\mathcal D}]$,  becomes a module over the 
 Hecke algebra  ${\mathcal H}={\mathbb Z}[q, q^{-1}]\otimes_{\mathbb Z} {\mathbb Z}[W]$, that is, the set 
 ${\mathcal H}$ of formal linear combinations of elements in $W$  with coefficients in 
 ${\mathbb Z}[q, q^{-1}]$.  If we denote by $T_w$ an element $w$ viewed inside  ${\mathcal H}$,  then the action of  $T_w$ on 
 any local system is determined by formulas describing the action of the  $T_{s_i} {\mathcal L}$ for the simple roots $s_i$. 
 These formulas are given  in p.371 of \cite{lusztig83} and if $w=s_{j_1}\cdots s_{j_k}$ is a reduced 
 expression and $m\in {\mathcal M}$, then $T_w m= T_{s_{j_1}} \cdots T_{s_{j_k}} m$ is independent of 
 the reduced expression.  
  
 We now express  $z(w,\mathcal L)$ in terms of the action of the Hecke algebra operators $T_w$.   If ${\mathcal L}\in {\mathcal D}$ is a local system on ${\mathcal O}_o$,  we have that  $T^{-1}_{w^{-1}} {\mathcal L} $ is a linear combination of  elements in the  set of    ${\mathcal D}$. Under our assumptions (i.e. real split Lie algebras and groups)   $T^{-1}_{w^{-1}} {\mathcal L} $, which gives the graded character of a principal series module relative to its weight filtration, contains exactly one  $K({\mathbb C} )$-equivariant  local systems on the open orbit ${\mathcal O}_o$ with non-zero coefficient (see \cite {casian86,casian99}).  This coefficient is $\pm q^{-z(w, {\mathcal L})}$.  This implies  that  $z(w,\mathcal L)$ is independent of the choice of a reduced expression (because $T_w$ is independent of a reduced expression). 
 \end{Remark}
 \begin{Example}\label{zz1}
In the case of $SL(2;{\mathbb R})$ (see Example \ref{list1}),  
 ${\mathcal D}= \{ {\mathcal C},{\mathcal L},  \delta_+, \delta_-  \}$ where  ${\mathcal C}$  denotes a trivial 
 sheaf  on ${\mathcal O}_o={\mathbb C}^*$, $\delta_{\pm }$ are sheaves supported on  the points $0$, 
 $\infty$ respectively,  and ${\mathcal L}$ is a non-trivial local system on ${\mathbb C}^*$. We have  
 $T_{s_1} {\mathcal C} = (q-2){\mathcal C} +(q-1)(\delta_-+\delta_+)$ and $T_{s_1}{\mathcal L}=-{\mathcal L}$. We now compute the numbers $z(s_1,\cdot)$ in terms of the Hecke algebra operators. In the Hecke algebra $T_{s_1}^{-1}= q^{-1}(T_{s_1} + (1-q))$. Hence by  applying 
  $T^{-1}_{s_1}$ to ${\mathcal C}$ we obtain  $-q^{-1} {\mathcal C} + q^{-1}(q-1)(\delta_++\delta_-)$. Thus the  coefficient of ${\mathcal C}$ is  $- q^{-z(s_1, {\mathcal C})}$ with $z(s_1, {\mathcal C})=1$.  We also have 
 $ q^{-1}(T_{s_1}{\mathcal L} + (1-q){\mathcal L})= q^{-1}(-{\mathcal L}+ (1-q){\mathcal L})=-{\mathcal L}$. Hence  $z(s_1, {\mathcal L})=0$. Note that these numbers extracted from the Hecke algebra operators defined in \cite{lusztig83}   correspond to the number of blow-ups in the  Toda lattice  in the $A_1$ case as  represented in Figure \ref{A1:fig} (in $\Gamma_{-}$, $\Gamma_{+}$ respectively).  See  Corollary \ref{zz2} below for a general statement.
 
 Finally, we note the connection with the ${\mathfrak sl}(2,{\mathbb R})$ 
 example in Introduction.  By rewriting   $T^{-1}_{s_1}q^1{\mathcal C}$ as
  $-q^{-1/2}\hat C_{\mathcal C} +\hat C_{\delta_{+}} +\hat C_{\delta_{-}}$ with $\hat C_{\mathcal C}=q^{-1/2}( {\mathcal C}+\delta_{+}+\delta_{-})$,  $\hat C_{\delta_{\pm}}=\delta_{\pm}$,  we recover the weight filtration of the module denoted $F^{even}(S^1)$ (replace $\hat C_{\mathcal C}$ with $C$ and $\hat C_{\delta_{\pm}}$ with $D_{\pm}$ and shift the weights there by $q^{-1/2}$).  The structure of the module $\Omega^{even}$ arises from  ${\mathcal C}$, and  the (irreducible) $F^{odd}(S^1)$ corresponds to ${\mathcal L}$.
 
 \end{Example}

 \begin{Proposition} \label{indepz} The numbers  $z(w,\mathcal L)$ are independent of the reduced expression used in Definition \ref {zw}.
 \end{Proposition}
 
 \begin{Proof}
 The argument follows from the formula in  Remark \ref{zz} above expressing $z(w,\mathcal L)$ in terms of the action of a Hecke algebra operator $T_w$, action which is independent of a reduced expression.  The number $q^{z(w,\mathcal L)}$ with $q$ a power of a prime is also  interpreted below  (and in \cite {casian99}) as a  Frobenius eigenvalue attached  to a Bruhat cell, that is, attached  to a Weyl group element $w$. Again from this second  interpretation, it follows that  the number $z(w,{\mathcal L})$  is  independent  of a reduced expression.  These two arguments are  related because the formulas  
 in p.371 of \cite{lusztig83} ultimately  arise from calculations that involve Frobenius eigenvalues.
   \end{Proof}

  
 We can now rewrite the description of the graph of incidence numbers  in \cite{casian99} as follows:

 \begin{Theorem} \label{casian99theorem} 
 The graph of incidence numbers ${\mathcal G}_{\mathcal L,B}$
 coincides with the graph obtained from Definition \ref{graphdef}.
  \end{Theorem}
 \begin{Proof}
 This theorem is proved in \cite{casian99} in representation theoretic terms and
corresponds to a similar result in \cite{koch} written in terms of different parameters. 
We give a geometric argument only  in the case $w\Rightarrow ws_i$.
Consider ${\mathcal B}_w \cup {\mathcal B}_{ws_i }$ a union of two $N$-orbits which
can be seen as a circle bundle over  ${\mathcal B}_w $. We fix a circle, a fiber of $\pi_i$  containing $x_w$,
and regard this set ${\mathcal B}_w \cup {\mathcal B}_{ws_i}$, this time 
as a  bundle over this circle. The fibers correspond to copies of 
${\mathcal B}_w $. By considering the volume
forms $\Lambda^{l(w)}w{\mathfrak n}\cap {\mathfrak n}^*$ the bundle becomes a line 
bundle over a circle, that is, either a cylinder or a M\"obius band. Determining the coefficient of the boundary of  ${\mathcal B}_{ws_i} $ along  ${\mathcal B}_w$ then becomes a question about
the triviality or non-triviality of this bundle. If the bundle is non-trivial
then there is a change in the orientation of the volume form around the circle
and this gives rise to a $\pm 2$.  
 \end{Proof}

\begin{Example}\label{exampleincidence} Consider $G=SL(3;{\mathbb R})$, $K=SO(3)$ and ${\mathcal L}$ as determined
 by the character defined by the action of
$T$ acting by $e_1e_2^{-1}$ (action on root vector for $e_1-e_2$). Since we are applying the character
to $e_i=\pm 1$, we can just write this character as $\chi ({\mathcal L})= e_1e_2$. We then have $s_1(e_1e_2)=e_1e_2$, $s_2(e_1e_2)=e_1e_3$, $s_2s_1(e_1e_2)=e_1e_3$, $s_1s_2(e_1e_2)=e_2e_3$ and
$s_1s_2s_1 (e_1e_2)=e_2e_3$. We write the characters of $\Lambda^{l(w)}{\mathfrak n}\cap {w^{-1}\mathfrak n}^*$, and obtain, for $w=e$, ${\mathfrak n}\cap {\mathfrak n}^*=0$ with $\chi_{e}=1$;
for $w=s_1$, ${\mathfrak n}\cap {s_1\mathfrak n}^*={\mathbb R} X_{e_1-e_2}$ with $\chi_{s_1}=e_1e_2$; for 
$w=s_2$, $ {\mathfrak n}\cap {s_2\mathfrak n}^*={\mathbb R} X_{e_2-e_3}$ with $\chi_{s_2}=e_2e_3$;
for $w=s_1s_2$, $ {\mathfrak n}\cap {s_2s_1\mathfrak n}^*={\mathbb R}X_{e_2-e_3}\oplus {\mathbb R}X_{e_1-e_3}$ with $\chi_{s_1s_2}=e_1e_2$; for
$w=s_2s_1$, $ {\mathfrak n}\cap {s_1s_2\mathfrak n}^*={\mathbb R}X_{e_1-e_2}\oplus{\mathbb R}X_{e_1-e_3}$ with $ \chi_{s_2s_1}=e_2e_3$; and for
$w=s_1s_2s_1$, $ {\mathfrak n}\cap {s_1s_2s_1\mathfrak n}^*={\mathbb R}X_{e_1-e_2}\oplus{\mathbb R}X_{e_2-e_3}\oplus{\mathbb R}X_{e_1-e_3}$ with $\chi_{s_1s_2s_1}=1$. Evaluating  $\Phi_{\alpha_1} ( \left(
\begin{matrix}
-1 & 0   \\
0 & -1 \\
\end{matrix}\right))$ in these characters  then means substituting $e_1=-1, e_2=-1, e_3=1$  and evaluating  $\Phi_{\alpha_2} ( \left(
\begin{matrix}
-1 & 0   \\
0 & -1 \\
\end{matrix}\right))$ corresponds to substituting $e_1=1, e_2=-1, e_3=-1$. For example 
$\epsilon (s_1s_2)$ is the pair $(-e_1e_2, -e_1e_2)$ in which we replace $e_1=-1,e_2=-1$
in the first coordinate and $e_1=1, e_2=-1$ in the second thus obtaining $(-+)$.  Therefore we obtain: 
$\epsilon (e)=(--)$, $\epsilon (s_1)=(-+)$,
$\epsilon (s_2)=(+-)$, $\epsilon (s_1s_2)=(-+)$,
$\epsilon (s_2s_1)=(+-)$, $\epsilon (s_1s_2s_1)=(--)$. From here it follows that the graph ${\mathcal G}_{\mathcal B}$  contains the
following edges, $s_1 {\Rightarrow} s_1s_2$ and $s_2{\Rightarrow} s_2s_1$. Therefore the graph ${\mathcal G}_{\mathcal B}$  which is the one derived earlier from the blow-ups of the Toda lattice in  Example \ref{exam}, i.e. ${\mathcal G}_{\mathcal B}={\mathcal G}$.

We now  tensor with  $w^{-1}\chi({\mathcal L})$ and obtain the list of
 $\chi_w({\mathcal L}) =w^{-1}\chi({\mathcal L}) \otimes \chi_w$:
$\chi_e ({\mathcal L})=e_1e_2$,
$\chi_{s_1} ({\mathcal L})=e_1e_2 \otimes e_1e_2=1 $,
$\chi_{s_2} ({\mathcal L})=e_2e_3 \otimes e_1e_3=e_1e_2 $,
$\chi_{s_1s_2} ({\mathcal L})=e_1e_2 \otimes e_1e_3=e_1e_2=e_2e_3 $,
$\chi_{s_2s_1} ({\mathcal L})=e_2e_3 \otimes e_2e_3=1 $.
$\chi_{s_1s_2s_1} ({\mathcal L})=1 \otimes e_2e_3=e_2e_3 $. Therefore we obtain: 
$\epsilon(e,{\mathcal L})=(-+)$,
$\epsilon(s_1,{\mathcal L})=(--)$,
$\epsilon(s_2,{\mathcal L})=(-+)$,
$\epsilon(s_1s_2,{\mathcal L})=(+-)$,
$\epsilon(s_2s_1,{\mathcal L})=(--)$,
$\epsilon(s_1s_2s_1,{\mathcal L})=(+-)$.

Now the graph ${\mathcal G}_{\mathcal L, B}$ has the following edge:
$e \overset{\mathcal L}{\Longrightarrow} s_2$. Since $s_1$ increases the length of both vertices
we also must have
$s_1 \overset{\mathcal L}{\Longrightarrow} s_2s_1$ and using that
 $s_2$ increases the length of both vertices 
 $s_1s_2 \overset{\mathcal L}{\Longrightarrow} s_2s_1s_2$.
We obtain the second  graph ${\mathcal G}_{(-+)}$ discussed in Example \ref{exam}.
\end{Example}

\section{Proof of Theorem \ref{mainT}}\label{mainTproof}
Now we give a proof of Theorem \ref{mainT} which provides the connection between the Toda lattice and the cohomology
of real flag manifold by showing that the graphs ${\mathcal G}_{\epsilon}$
are the graphs of the incidence numbers for cohomology of the flag manifold
$\check{\mathcal B}$. Let us begin with the following Proposition:

\begin{Proposition}\label{propsigns} For  $\epsilon(w, {\mathcal L})=(\epsilon_1\ldots \epsilon_l)$,
we have $\epsilon(ws_i,{\mathcal L})=(\epsilon'_1\ldots \epsilon'_l)$ where
$\epsilon'_j =\epsilon_j\epsilon_i^{C_{i,j}}$. 
\end{Proposition}
\begin{Proof}
We have 
${\mathfrak n}\cap {s_iw^{-1}\mathfrak n}^*={\mathfrak n}\cap {w^{-1}\mathfrak n}^*\oplus {\mathbb R} X_{\alpha_i}$,
where $X_{\alpha_i}$ is a root vector.  Therefore
$\Lambda^{l(w)+1}{\mathfrak n}\cap {s_iw^{-1}\mathfrak n}^*=\Lambda ^{l(w)}{\mathfrak n}\cap {w^{-1}\mathfrak n}^*\wedge  {\mathbb R} X_{\alpha_i}$.
Hence the character of $T$ given by $\Lambda ^{l(w)}{\mathfrak n}\cap {w^{-1}\mathfrak n}^*$ changes exactly
by the character $\chi_{\alpha_i}$  corresponding to the action of $T$ on the root vector $X_{\alpha_i}$.
The evaluation of  $\chi_{\alpha_i}$ on  $\Phi_{\alpha_j} ( \left(
\begin{matrix}
-1 & 0   \\
0 & -1 \\
\end{matrix}\right))$ is given by $e^{<\alpha_i; \sqrt{-1}\pi h_{\alpha_j }   >}$
that is as $e^{\sqrt{-1}\pi C_{i,j}}$. Hence we obtain $\epsilon_j \to \epsilon_j\epsilon_i^{C_{i,j}}$.
\end{Proof}

Note that the action of $s_i$ on the signs in Proposition \ref{propsigns} is different from the action of $s_i$ on signs of the Toda lattice in Definition \ref{act}. The difference is that in the case of the Toda lattice with $C_{j,i}$, the transpose of $C_{i,j}$ is involved; because of this, the Lie algebra $\check {\mathfrak g}$ is involved.    
Thus the set of signs $\{  \epsilon(w, {\mathcal L})\,|\,w\in W \}$ is the $W$-orbit of the signs $\epsilon(e, {\mathcal L})$ with the $W$-action defined with respect to  the Lie algebra $\check {\mathfrak g}$ $({\rm Definition}\, \ref{act})$.  This is also
the set of signs $({\rm sgn}(a_1)\ldots{\rm sgn}(a_l))$ that occur  in the polytope $\Gamma_{\epsilon  (e, {\mathcal L})}$  with respect to the  Lie algebra $\check {\mathfrak g}$.
For example, in the case in which ${\mathfrak g}$ is  a Lie algebra of type $B$, the signs in a polytope in the isospectral manifold of the Toda lattice correspond to a Lie algebra of type $C$.

The following proposition will then relate the Frobenius eigenvalues $q^{z(w, {\mathcal L})}$ appearing in the cohomology of ${\mathcal B}(k_q)$ with the numbers  $\check \eta (w, \epsilon  ({\mathcal L}, w) )$ for the Toda lattice associated to the Lie algebra $\check {\mathfrak g}$.

\begin{Proposition}\label{equalweights} Assume that ${\mathfrak g}$ is a semisimple real split Lie algebra  then $z(w, {\mathcal L})=\check \eta (w, \epsilon  (w,{\mathcal L}) )$.
\end{Proposition}
\begin{Proof} By Proposition \ref{propsigns} the set of signs  $ \{ \epsilon  (w,{\mathcal L})\,|\,w\in W \}$
keeping track of triviality or non-triviality of the local systems  ${\mathcal L}_w$ along certain
fibers of $\pi_i$ is just the $W$-orbit of the sign $\epsilon  (e, {\mathcal L})$. This is just
the set of signs $({\rm sgn}(a_1)\ldots {\rm sgn} (a_l))$ in the polytope $\Gamma_{\epsilon  (e, {\mathcal L})}$
for the Lie algebra $\check {\mathfrak g}$  as noted above.   
 We then note that the condition for $w\Rightarrow ws_i$ in Definition \ref{eta}  that $\epsilon_i=+$
corresponds to the non-triviality of the sheaf ${\mathcal L}_w$ along the fiber of $\pi_i$
containing $x_w$. From here it follows that $z(w, {\mathcal L})=\check \eta (w, \epsilon  (w, {\mathcal L}) )$.
\end{Proof}

We obtain the following expression for  the numbers $\check \eta(w,\epsilon)$ in terms of Hecke algebra operators:

\begin{Corollary} \label{zz2}  If  $\epsilon=\epsilon(e,{\mathcal L})$, then  $\pm q^{-\check \eta(w,\epsilon)}$ is the coefficient of ${\mathcal L}$ in $T_{w^{-1}}^{-1}{\mathcal L}$ in the module over the Hecke algebra in \cite{lusztig83} (Remark \ref{zz}).
\end{Corollary}
\begin{Proof} This follows from Remark  \ref{zz} coupled with Proposition \ref{equalweights}.
\end{Proof}

In particular we recover the number of blow-ups  $\eta(s_1,-)=\eta(s_1)=1$ and  $\eta(s_1,+)=0$ in the ${\mathfrak sl}(2;{\mathbb R})$  case from the action of the Hecke algebra operator $T_{s_1}$ as in Example \ref{zz1}.

\begin{Corollary}\label{indep} The number $\eta(w,\epsilon)$ is independent of a reduced expression of $w$. 
\end{Corollary}
\begin{Proof} This follows  from Corollary \ref{zz2}  or directly from Proposition \ref {equalweights} because the numbers  $z(w, {\mathcal L})$ are independent of a reduced expression of $w$ (Proposition \ref{indepz}).  It is possible to choose the group $G$ so that  all possible $\epsilon$ are of the form  $\epsilon  (e, {\mathcal L})$ ($\tilde G$ in Remark \ref{group1}  or in  \cite {casian:02}).  
\end{Proof}

The following Proposition completes the proof of Theorem \ref {mainT} (see Definition \ref{graph}
for $\mathcal G_{\epsilon}$):

\begin{Proposition}\label{graphToda} We have  $w_1\, { \overset{\mathcal L}{\Longrightarrow }} \,w_2$ if and only if 
\begin{itemize}
\item[a)]   $w_1\le w_2$ in the Bruhat order,
\item[b)]  $ l(w_2)=l(w_1)+1$\,,
\item[c)]  $\check \eta (w_1, \epsilon (e, {\mathcal L}))=\check \eta (w_2, \epsilon (e, {\mathcal L}))$\,,\\
\item[d)] $\epsilon(w_1,{\mathcal L})=\epsilon(w_2,{\mathcal L})$\,.
\end{itemize}
\end{Proposition}
\begin{Proof} Note that  $w_1 \le w_2$ in the Bruhat order  with  $l(w_2)=l(w_1)+1$ if and only if there are $a,x, s_i\in W$ such that $w_1=ax$, $w_2=as_ix$ and $l(axs_i)=l(a)+l(x)+1$.   Hence, by definition, if a) and b) are satisfied  $w_1 \overset {\mathcal L}{ \Longrightarrow } w_2$ if and only if  $a \overset{\mathcal L}{\Longrightarrow } as_i$.  In turn this is true if and only if  $z(a , \epsilon  (a,{\mathcal L}) )=z(as_i , \epsilon  (a,{\mathcal L}) )$.  By Proposition  \ref{equalweights}, this is true  if and only if  $\check \eta (a, \epsilon (e, {\mathcal L}))=\check \eta (as_i, \epsilon (e, {\mathcal L}))$, and moreover, if and only if  the signs  $\epsilon  (a, {\mathcal L})$ and $\epsilon  (as_i, {\mathcal L})$ agree and the $i$th sign is $+$ on both sides. Again by Proposition \ref{equalweights}, since $\epsilon  (ax, {\mathcal L})$ and $\epsilon  (as_ix, {\mathcal L})$ are obtained by applying $x^{-1}$ to $\epsilon  (ax, {\mathcal L})= \epsilon  (as_ix, {\mathcal L})$, using the $W$-action we have $\check \eta (a, \epsilon (e, {\mathcal L}))=\check \eta (as_i, \epsilon (e, {\mathcal L}))$ if and only if   $\check \eta (ax, \epsilon (e, {\mathcal L}))=\check \eta (as_ix, \epsilon (e, {\mathcal L}))$. Thus we conclude that whenever a) and b) are satisfied  $w_1 \overset{\mathcal L}{\Longrightarrow } w_2$ if and only if  $\check \eta (w_1, \epsilon (e, {\mathcal L}))=\check \eta (w_2, \epsilon (e, {\mathcal L}))$.
\end{Proof}

\smallskip

\begin{Remark}
 Let $W^{-}=\{ w\in W\,|\, w^{-1}(- \ldots  -)=(-\ldots -) \}$. Note that  under  the correspondence between local $K$-equivariant local systems and signs ${\mathcal L}\to \epsilon (w, {\mathcal L})$, a trivial sheaf ${\mathcal L}$ corresponds to an element in $W^{-}$.    By Theorem \ref{mainT}  the vertices of the graph of incidence numbers belonging to $W^{-}$ are the only ones that  may contribute to  rational cohomology.   
\end{Remark}


\section{Rational cohomology of $K$ and $p(q)$} 
Here we first show that the number of blow-ups $\eta(w)$ for each $w\in W$
is given by the eigenvalues of the Frobenius map on the cohomology of flag manifold.
We then show that the polynomial $p(q)$ defined in Definition \ref{pq} is related to the
order of Chevalley group $K({\mathbb F}_q)$.

\subsection{Frobenius eigenvalues and $\eta(w)$}
Let us consider a filtration,
\[ 
\emptyset\, \subset \, {\mathcal Y}_0\, \subset  \, {\mathcal Y}_{1} \, \subset \, \cdots\, \subset \, {\mathcal Y}_{l(w_*)}= {\mathcal B}_{\mathbb C}
\] given by intersection of ${\mathcal O}_o$
 with  the $N({\mathbb C})$ cells
inside $G({\mathbb C})/B({\mathbb C})$.  We obtain  coboundary maps,
$H^j({\mathcal Y}_j, {\mathcal Y}_{j-1}; {\mathbb C}) \to 
H^{j+1}({\mathcal Y}_{j+1}, {\mathcal Y}_{j}; {\mathbb C})$ which
give rise to a chain complex computing the cohomology of ${\mathcal O}_o$.
For example in the case of $SU(1,1)$ above,  ${\mathcal Y}_0=\{ \infty \}$ and ${\mathcal Y}_1={\mathbb C} \setminus \{ 0 \}$. 
Each $w$ corresponds to a dual of a Bruhat cell and  contributes to the cohomology of  $H^{l(w)}({\mathcal Y}_{l(w)}, {\mathcal Y}_{l(w)-1}; {\mathbb C})$
giving rise to a cohomology class $[w]_{\mathbb C}$.
This can be done with etale cohomology with coefficients in $\bar {\mathbb Q}_m$
and a field of positive characteristic.  In this case a Frobenius action arises
in cohomology.
 
We have the following proposition  for the Frobenius eigenvalue of the cohomology
class $[w]_{k_q}$ in etale cohomology over a field $k_q$ algebraic closure of ${\mathbb F}_q$ of characteristic $p$.  This assumes $\bar {\mathbb Q}_m$ coefficients where $m$ is relatively prime to $p$ and $p\not=2$ and $x^2+1$ factors over ${\mathbb F}_q$.

\begin{Proposition}\label{frob}
The cohomology class $[w]_{k_q}$ in $H^{l(w)} ({\mathcal Y}_{l(w)},  {\mathcal Y}_{l(w)-1}; \bar {\mathbb Q}_m)$ corresponding to $w\in W$
has Frobenius eigenvalue given by  $q^{\check \eta(w)}$. 
\end{Proposition}
\begin{Proof}
This statement is obtained by expressing Proposition 9.5  in \cite{casian99} in terms
of new notation where local systems are expressed in terms of their signs $\epsilon (w,{\mathcal L})$.  This corresponds to Corollary \ref{zz2}. The argument  in Proposition 9.5  in \cite{casian99} for the real split case reduces to the ${\mathfrak sl}(2;{\mathbb R})$ case and, strictly speaking, gives Frobenius eigenvalues  up to a sign $\pm$. This choice of a sign corresponds to the two possibilities $q-1$ or $q+1$ for the number of points in $SO(2;{\mathbb F}_q)$. By assuming that ${\mathbb F}_q$ contains $\sqrt{-1}$, the formula for the number of points is fixed as $q-1$, and thus the sign is fixed too as in Example \ref {sl2x} below, or the $A_1$ example given in  Introduction .

The upshot of this is that  the number of blow-ups $\eta(w)$ associated to a vertex $w$ in
the Toda lattice and the Frobenius eigenvalue of
 $[w]_{k_q}$ in $H^{l(w)}({\mathcal Y}(k_q)_{l(w)}, {\mathcal Y}(k_q)_{l(w)-1}; \bar {\mathbb Q}_m)$
 are given by the same formula in Definition \ref{eta}. 
 \end{Proof}

\begin{Example}\label{sl2} The case of $SU(1,1)$: Recall in Example \ref{suoneone}
that $SU(1,1)$ corresponds to the real form of $SL(2;{\mathbb C})$. In this case the number of ${\mathbb F_q}$ points in
${\mathcal O}_o(k)$ equals $q-1$. This corresponds to the following Frobenius eigenvalues
$1$ on  $H_c^1(k_q^*;\bar {\mathbb Q}_m)$ and $q$ on $H_c^2(k_q^*;\bar {\mathbb Q}_m)$.
The map $Fr$ is given by $Fr: z\mapsto z^q$.
\end{Example}
\begin{Example}\label{sl2x} The case of $SL(2;{\mathbb R})$:
The real form and  $K$ are determined by 
 $\theta (A)= A^*$, the transpose of the inverse of $A$. Hence
$K(k)$ is given by  $\left(
\begin{matrix}
x & -y   \\
y &  x \\
\end{matrix}\right) $ where $x^2+y^2=1$.  Roots of $x^2+1$ are necessary to diagonalize $K$
and this introduces a difference between fields ${\mathbb F}_q$ such that $x^2+1$
splits into linear factors  and fields for which this is not the case. We have an
action of  $SL(2;k_q)$  on ${\mathbb P}^1$ which can be described with fractional linear
transformations. Hence $K$ acts as follows on $z$,
$\left(
\begin{matrix}
x & -y   \\
y &  x \\
\end{matrix}\right)\cdot z=\frac{xz-y}{yz+x}$.
When $z=0$ we obtain the orbit $\{r= -\frac{y}{x} \,|\,x^2+y^2=1, x\not=0 \} \cup \{\infty \}$.
We have $r=\pm \frac{y}{\sqrt{1-y^2}}$ and $y=\pm \frac{r}{\sqrt {1+r^2}}$ in the algebraically closed field $k_q$.  Therefore for any $r$ such that $r^2+1\not=0$, we get a $y$ and then an $x=\pm \sqrt{1-y^2}$.
Hence the set of $r\in k_q$, the points in the orbit ${\mathcal O}_o(k_q)$ which are not $\infty$,
is given as the set $\{ r\,|\,  r^2+1\not=0 \}$.  We now compute the number of ${\mathbb F}_q$ points (see the $A_1$ example in Introduction).  
There are two possibilities depending on the 
finite field ${\mathbb F}_q$: if $r^2+1\not=0$ in ${\mathbb F}_q$ we obtain a contribution of  $q$ elements corresponding to all the possible values of $r$, and if we include the point $\infty$, we find
 $q+1$ elements in  ${\mathcal O}_o({\mathbb F}_q)$. In the case when
$r^2+1=0$ has two solutions, $\pm \sqrt {-1}$, we obtain $q-2$ points from the possible values of $r$
excluding these two roots, then taking into account the point  $\infty$ we have a total of  $q-1$ points. These formulas corresponds to the following Frobenius eigenvalues,
$\pm 1$ on  $H_c^1({\mathcal O}_o(k);\bar {\mathbb Q}_m)$ and $q$ on $H_c^2({\mathcal O}_o(k_q);\bar {\mathbb Q}_m)$. The eigenvalue $+1$ of Frobenius acting on $H_c^1({\mathcal O}_o(k);\bar {\mathbb Q}_m)$ then corresponds to the case when $x^2+1$ factors over ${\mathbb F}_q$.
By considering ${\mathbb F}_p[\sqrt{-1}]$ when necessary, it is possible to assume that  $x^2+1$ splits
in ${\mathbb F}_q$.
For example over ${\mathbb F}_3$ ($q=3)$ the equation $x^2+1=0$ has no roots and this is
the reason that $x_1^2+x_2^2=1$ contains $4=q+1$ points, $\{ (1,0),(2,0),(0,1),(0,2) \}$. However we can change  this  by simply considering  ${\mathbb F}_{3^2}={\mathbb F}_3[\sqrt{-1}]$ and  then obtain $q-1=8$ points for the new $q=3^2$. Over
${\mathbb F}_5$, the equation $x^2+1=0$ has two solutions $ x=2, 3$. In this case $x_1^2+x_2^2=1$
contains $4=q-1$ points, $\{ (1,0),(4,0),(0,1),(0,4) \}$. 
\end{Example}


 By using the Lefschetz Fixed Point  Theorem applied to the 
 Frobenius map to count
 the number of ${\mathbb F}_q$ points  and that the number $|K({\mathbb F}_q)|$ is given in p.75 of  \cite{carter}, we  obtain the following:

\begin{Proposition}\label{pointsinflag} For any $q$, power of a prime $p\not=2$ such that $x^2+1$ factors  as product of linear terms, we have  $|\check { \mathcal O}_o({\mathbb F}_q)|=q^{r_1}p(q)$ with $r_1={{\rm dim}(\check K)- {\rm deg} (p(q))}$ where $p(q)$ is given by
\[
p(q)=(-1)^{l(w_*)}\sum_{w\in W^-} (-1)^{l(w)} q^{ \eta(w)}\,.
\]
Here $W^-$ is defined by $W^{-}=\{ w\in W\,|\, w^{-1}(- \ldots  -)=(-\ldots -) \}$.
\end{Proposition}
\begin{Proof} The cohomology of $\check {\mathcal O}_o$ is given in terms
of a chain complex 
$$\cdots \to H^s(\check {\mathcal Y}_s, \check {\mathcal Y}_{s-1};{\mathbb C}) \to 
H^{s+1}(\check {\mathcal Y}_{s+1}, \check {\mathcal Y}_{s};{\mathbb C}) \to \cdots $$
Take ${\mathcal L}$ the trivial local system. Those $w$ for which ${\mathcal L}_w$ 
is trivial correspond to $w\in W^-$.  These are the $w$ which contribute
to the Betti numbers. We now recall that from the Lefschetz fixed point formula for $Fr$ the alternating sum is given by
\[
\sum_{s}(-1)^s {\rm Tr}\left((Fr)_*|_{H_c^s(\check {\mathcal O}_o(k_q);\bar {\mathbb Q}_m)}\right)=|\check {\mathcal O}_o(k_q)|
\]
By the Poincare duality, $H_c^{2{\rm dim}( \check K)-s}(\check {\mathcal O}_o(k_q); \bar {\mathbb Q}_m)$ is the dual of
$H^{s}(\check {\mathcal O}_o(k_q); \bar {\mathbb Q}_m)$. 
Then the alternating sum of the traces of Frobenius in cohomology is 
$(-1)^{l(w_*)}p(q)=\sum_{w\in W} (-1)^{l(w)} q^{\eta (w)}$. Since only elements in $W^-$ contribute to
the Betti numbers, the $W$ in the sum can be replaced by $W^-$. Replacing cohomology with cohomology with proper supports amounts to replacing $q\to q^{-1}$ and multiplying by $q^{{\rm dim} (\check K)}$ this previous expression.  We write ${\rm dim} (\check K)= r_1+ \eta(w_*)$ and now 
$\sum_{w\in W} (-1)^{l(w)} q^{\eta (w)}$ gives rise to 
$q^{r_1}\sum_{w\in W} (-1)^{l(w)} q^{\eta (w_*)-\eta(w)}$. 
By Proposition \ref{poincare}  ($\eta(w_*w)=\eta(w_*)-\eta(w)$),
we obtain  $q^{r_1}\sum_{w\in W} (-1)^{l(w)} q^{\eta (w_*w)}$. Since 
$l(w_*w)= l(w_*)-l(w)$, our expression becomes 
$q^{r_1}(-1)^{l(w_*)} \sum_{w\in W} (-1)^{l(w_*w)} q^{\eta (w_*w)}$.
This is $q^{r_1}p(q)$. Therefore the alternating sum of the traces of $Fr$ in cohomology with
proper supports is $q^{r_1}p(q)$. Using the Lefschetz fixed point formula for $Fr$ we conclude
$|\check {\mathcal O}_o(k_q)|=q^{r_1}p(q)$.
\end{Proof}
\begin{Notation}
We now consider several projections:
\begin{itemize}
\item{}
We first express $p: G({\mathbb C})/B({\mathbb C})\to \{* \}$, the projection to a point $\{ *\}$ as a composition of the projection, 
 $\pi_i: G({\mathbb C})/B({\mathbb C}) \to G({\mathbb C})/P_i({\mathbb C})$ and the projection to a point
$p_i: G({\mathbb C})/P_i({\mathbb C})\to \{ * \}$ so that $p=p_i\circ \pi_i$. 
\item{}
Denote  $p^o$ the restriction of $p$ to the $K({\mathbb C})$-orbit  ${\mathcal O}_o$.
Similarly denote $p^o_i$ the restriction to ${\pi_i({\mathcal O}_o)}$, that is,  $p_i^o=p_i|_{\pi_i({\mathcal O}_o)}$ and finally, $\pi^o_i=\pi_i|_{{\mathcal O}_o}$ the restriction of $\pi_i$ to  ${\mathcal O}_o$.
\item{}
The fibration $\pi_i$ has fibers $E_x$ identifiable with ${\mathbb P}^1$ and the fibration $\pi_i^o$ has
fibers $E_x^o$ identifiable with ${\mathbb C}^*= {\mathbb C}\setminus \{ 0 \}$. 
\end{itemize}
\end{Notation}

In the proof of the following proposition we use a sheaf-push-forward of the form $f_{!}$ (see chapter VII of \cite{iversen}). Recall that this notation refers to  a push-forward with compact supports. Moreover we consider derived functors (p.51 of \cite{iversen}), obtained by considering injective resolutions (or soft resolutions) and then applying the push-forward to it.  The reason for using this construction here is that the ordinary  cohomology of an a space $X$ is  obtained by applying a push-forward  $f: X\to \{ * \}$ (to a point),  to a constant sheaf.  More precisely, we consider the derived functors $R^sf_*$  and apply them to a constant sheaf.  If the push-forward is with proper supports,  $R^sf_!$, then one obtains cohomology with compact supports. These two notions of cohomology are related by Poincare duality when $X$ is a smooth algebraic variety (Theorem 2.6 in p.263 of \cite{iversen}).

\begin{Proposition}\label{flagvsgroup} We have $H^k(K; {\mathbb Q})= H^k(K/T; {\mathbb Q})$ for all k.
Equivalently, for any nontrivial $K({\mathbb C})$-equivariant local system ${\mathcal L}$, we have
$H^*({\mathcal O}_o, {\mathcal L} )=0$.
\end{Proposition}
\begin{Proof} The fibration $s:K({\mathbb C}) \to  K({\mathbb C})/ T$ has finite fibers
and  $s_*(1)$ is a direct sum of all $K({\mathbb C})$-equivariant local
systems on  $K({\mathbb C})/ T$. Hence the rational cohomology of $K$  is a direct sum of 
cohomology groups of $K/T$ with twisted coefficients.   To  verify that $H^*(K/T;{\mathbb Q})=H^*(K; {\mathbb Q})$ then it suffices to  show that rational cohomology of $K/T$ with twisted coefficients ${\mathcal L}$ is zero. 

The proof now becomes  parallel,  to the proof  that  $p_\epsilon (q)=0$ when $\epsilon$ contains a positive sign $\epsilon_i=+$.   In the proof of Proposition \ref{pzero}  a single arrow $e\Rightarrow s_i$ produces  $p_\epsilon (q)=0$. Here the non-triviality of a local system ${\mathcal L}$ corresponds to $\epsilon_i=+$ for  the sign $\epsilon=\epsilon (e,{\mathcal L})$.

We have $R^sp^o_!({\mathcal L})=H_c^s({\mathcal O}_o; {\mathcal L})$, the cohomology
with proper supports and with local coefficients,  which can be written as  $H_c^*(K/T; {\mathcal L}$).
Similarly with the restriction $p^o_i$,
the composition of functors $p_i^o\circ \pi^o_i=p^o$
gives rise to the spectral sequence (composition of functors spectral sequence as in p.343 of \cite{cartan} or Theorem 7.15 in \cite{iversen} in terms of derived categories) $R^k(p_i^o)_!(R^s(\pi^o_i)_!{\mathcal L})\Rightarrow R^{k+s}p_!{\mathcal L}.$
To show that $R^{r}p_!{\mathcal L}=0$ for all $r$,  it suffices to show that $R^s(\pi^o_i)_!{\mathcal L}=0$
for at least one $i$ and all $s$, $s=1,2$.
For this  it  will be  enough to compute the stalks over one single  point $x$,  of the sheaves  $R^s(\pi^o_i)_!{\mathcal L}$. On the other hand,  {\it proper base change } (section III. 6 of \cite{iversen}) implies that  the stalks   $(R^{s}(\pi^o)_!{\mathcal L})_x$ are given as $H_c^{s}(E^o_x; {\mathcal L}|_{E^o_x})$, that is,  the cohomology with proper supports of 
 the local system ${\mathcal L}$ along a copy of  ${\mathbb C}^*={\mathbb C}\setminus \{0 \}$.
 This has now reduced the argument to  a calculation in the flag manifold for $A_1$, 
since a fiber $E_x$ of $\pi_i$ can be identified with ${\mathbb P}^1$
and cohomology is along ${\mathbb C}^*$, a $K({\mathbb C})$-orbit  for $K=SO(3)$ or $SU(2)$. Denote $\pi_i|_{E_x}=\pi_x$. Then we must compute $R^*(\pi^o_x)_!{\mathcal L}|_{E^o_x}$. To show that this is zero for some $i$. It now  suffices that  ${\mathcal L}$ is non-trivial along the fiber $(\pi^o_i)^{-1}(x_e)$ for some $i$ (again,  corresponding to  $\epsilon_i=+$ with $\epsilon=\epsilon (e,{\mathcal L})$ ). Since the rational cohomology of such a  non-trivial local system along ${\mathbb C}^*$  is zero,  we  have that  $R^s(\pi^o_i)_!{\mathcal L}=0$ as needed in our proof. 
 \end{Proof}

\subsection{$|K({\mathbb F}_q)|$, $p(q)$ and the rational cohomology of $K$} \label{coh}

Recall that the polynomial  $|K({\mathbb F}_q)|$ has a factorization $|K({\mathbb F}_q)|=q^r \Pi_{i=1}^g (q^{d_i}-1)$  where $d_1,\ldots,d_g$ are degrees of basic Weyl group invariant polynomials for $K$
with $g={\rm rank}(K)$
and $r$ is the total number of positive roots associated with $K$ (see \cite{carter}). Therefore the polynomial  $|K({\mathbb F}_q)|$ encodes the information that describes the cohomology ring of $K$. This cohomology ring  is well-known; it is an exterior algebra with $g$ generators of degrees $2d_i-1$.   Furthermore, if we write the exterior algebra as $H^*(K;{\mathbb Q}) =\Lambda(q^{d_1}x_1,\ldots, q^{d_g}x_g)$, we obtain a filtration of cohomology  by the powers of $q$ that result. This filtration  can be seen to correspond to a filtration by Frobenius eigenvalues in etale cohomology of $K(k_q)$. We illustrate the relation between the polynomial $|K({\mathbb F}_q)|$ and cohomology with the examples
of $A_3$ and $G_2$:

\begin{Example}
   In the case of type $A_3$, we have $K=SO(4)$  and $|SO(4;{\mathbb F}_q)|=q^2(q^2-1)(q^2-1)$, i.e. $d_1=2, d_2=2$ (since $SO(4)$ and $SU(2)\times SU(2)$ have the same Lie algebra, and for $SU(2)$ we have $d_1=2$ associated to the Casimir operator).  We consider then the exterior algebra with generators $\{q^2x_1, q^2x_2\}$. From here one can write  the rational cohomology of $K$, and even more, $\bar {\mathbb Q}_m$ the etale cohomology of $K(k_q)$  (with $m\not=2$ a prime which is relatively prime to $q$),
\[
\left\{\begin{array}{llll}
H^0(SO(4;k_q) ; \bar {\mathbb Q}_m)=\bar {\mathbb Q}_m(q^01)\quad &:\quad {{\rm trace}=q^0}\\ 
H^3( SO(4;k_q);\bar {\mathbb Q}_m)=\bar {\mathbb Q}_m(q^2x_1)\oplus \bar {\mathbb Q}_m(q^2x_2) \quad &:\quad {{\rm trace}=2q^2}\\
H^6( SO(4;k_q) ;\bar {\mathbb Q}_m)=\bar {\mathbb Q}_m(q^4x_1 \wedge x_2) \quad &:\quad {{\rm trace}=q^4}
\end{array}\right.
\]

The alternating sum of the traces of Frobenius eigenvalues  in cohomology gives $p(q)=(q^2-1)^2$. If we consider cohomology with proper supports then we obtain $q^2(q^2-1)^2=|SO(4;{\mathbb F}_q)|$.  The generators of degree $3$ in this example  correspond to linear combinations  $x_1=[121]\pm [232]$, $x_2=[321]\pm [123]$ in terms of the graph in  Figure \ref{A3inc:fig} and the product $x_1\wedge x_2$ corresponds to the longest element  $w_*=[123121] $. The graph in Figure \ref{A3inc:fig}  is derived from the blow-ups of the Toda lattice but it  agrees with  the graph of incidence numbers of $G/B$  (by Theorem \ref{mainT} or  by  p.529 of  \cite{casian99}). These representatives have eigenvalues of Frobenius $q^2$, $q^2$ and $q^4$ respectively. 

 In the case of $G_2$ (Example \ref{A2G2}) we also have that the Lie algebra of $K$ agrees with the Lie algebra of $SU(2)\times SU(2)$. Therefore, in spite of the fact that the structure of the polytopes is very different  in $A_3$ and $G_2$ (i.e. different dimension, different number of vertices and different number of divisors),  the polynomial $|K({\mathbb F}_q)|=q^2(q^2-1)^2$.
\end{Example}

In the following, we assume that  $p\not=2$ and $x^2+1$ factors over the field $ {\mathbb F}_q$:

\begin{Theorem} 
\label{Kpq}
  The polynomial $p(q)$ satisfies  $ p(q)=\tilde p(q)$, which is $q^{-r}| \check K({\mathbb F}_q)|$ with
  $r={\rm dim} (\check K)- {\rm deg} (p(q))$. Moreover $p(q)$ factors as   $p(q)= \Pi_{i=1}^g (q^{d_i}-1)$
where $d_1,\ldots, d_g$ are degrees of basic Weyl group invariant polynomials for $\check K$
of rank $g$.  The polynomial $p(q)$  is given by the following explicit formulas:
 
  $A_{l}$:  $\check K=SO(l+1)$,
 \begin{itemize}
 \item[] $l$ even: $p(q)=(q^2-1)(q^4-1)\cdots (q^{l-2}-1)(q^{l}-1),\quad g=l/2$ 
 \item[] $l$ odd:  $p(q)=(q^2-1)(q^4-1) \cdots (q^{l-3}-1)(q^{l-1}-1)(q^g-1),\quad g=(l+1)/2$ 
 \end{itemize}

 $B_l$:  $\check K=U(l)$,
 \begin{itemize}
\item[]  $p(q)=(q-1)(q^2-1)(q^3-1) \cdots (q^l-1),\quad g=l$ \par
\end{itemize}

$C_l$:  $\check K=SO(l)\times SO(l+1)$,
 \begin{itemize}
\item[]  $l$ even:  $p(q)=(q^2-1)^2(q^4-1)^2 \cdots (q^{l-2}-1)^2(q^l-1)(q^{l/2}-1),\quad g=l$ \par
\item[]  $l$ odd:  $p(q)=(q^2-1)^2(q^4-1)^2 \cdots (q^{l-1}-1)^2(q^{(l+1)/2}-1),\quad g=l$ \par
\end{itemize}

$D_l$:  $\check K=SO(l)\times SO(l)$,
 \begin{itemize}
\item[]   $l$ even: $p(q)=(q^2-1)^2(q^4-1)^2(q^6-1)^2...(q^{l-2}-1)^2(q^{l/2}-1)^2,\quad g=l$ \par
\item[]   $l$ odd: $p(q)=(q^2-1)^2(q^4-1)^2(q^6-1)^2...(q^{l-1}-1)^2,\quad g=l-1$ \par
\end{itemize}
    
$E_6$:  ${\rm Lie}(\check K)=\mathfrak{sp}(4)$,
\begin{itemize}
\item[]   $p(q)=(q^2-1)(q^4-1)(q^6-1)(q^8-1),\quad g=4$
\end{itemize}

 $E_7$:   ${\rm Lie}(\check K)={\mathfrak{su}}(8)$,
\begin{itemize}
\item[]  $p(q)=(q^2-1)(q^3-1)(q^4-1)(q^5-1)(q^6-1)(q^7-1)(q^8-1),\quad g=7$ 
\end{itemize}

$E_8$: ${\rm Lie}(\check K)={\mathfrak{so}}(16)$,
\begin{itemize}
\item[]  $p(q)=(q^2-1)(q^4-1)(q^6-1)(q^8-1)(q^{10}-1)(q^{12}-1)(q^{14}-1)(q^8-1),\quad g=8$ \par
\end{itemize}

$F_4$:  ${\rm Lie} (\check K)={\mathfrak{sp}}(1)\times {\mathfrak{sp}}(3)$,
 \begin{itemize}
 \item[]  $p(q)=(q^2-1)^2(q^4-1)(q^6-1),\quad g=4$   
\end{itemize}

$G_2$:  ${\rm Lie}(\check K)=\mathfrak{su}(2)\times {\mathfrak{su}}(2)$,
  \begin{itemize}
\item[]  $p(q)=(q^2-1)^2,\quad g=2$ 
\end{itemize}

 \end{Theorem}
 
 \begin{Proof}  By Proposition \ref {pointsinflag}, we have $q^{r_1}p(q)=| \check {\mathcal O}_o({\mathbb F}_q)|$. Then by  Proposition \ref {flagvsgroup},  $| \check {\mathcal O}_o({\mathbb F}_q)|=| \check K({\mathbb F}_q)|$ (hence $r=r_1$).
 By possibly extending the field ${\mathbb F}_q$ as described in Remark \ref{so}, we have a factorization $|\check K({\mathbb F}_q)|=q^r\Pi_{i=1}^{g} (q^{d_i}-1)$ (see p.75 of  \cite{carter}). 
 Moreover the polynomials  $|K({\mathbb F}_q)|$ are also listed in  p.75 of  \cite{carter}. The only groups $K$ that occur for the real split semisimple Lie groups are such that by possibly replacing $q$ with $q^2$ as  in Remark \ref{so},  $K({\mathbb F}_q)$  is of type $A$, $B$, $C$ or $D$ in the classification of finite Chevalley groups  given in p.37 of \cite{carter}.   From the polynomials $|K({\mathbb F}_q)|$ one then  obtains the polynomials $p(q)$ by dividing by $q^{r}$. 
 Note also that where the group $K=U(n)$ appears, we are dealing with a Chevalley group of type $A$, i.e.  $SL(n; {\mathbb F}_q)$, after reduction to positive characteristic and then taking the ${\mathbb F}_q$ points. This is because $SU(n)$ has complexification $SL(n;{\mathbb C})$.  The remark applies to the case of a Lie algebra of type $B_l$ above and $E_7$ (e.g.  a maximal compact Lie subgroup for type $C_l$ is $U(l)$ and we need to use of  the check $\check{}$ which exchanges the maximal compact subgroups of $C$ and $B$).
 \end{Proof}

  From  here  we obtain: 

 \begin{Proposition} 
 \label{totalBnumber}
 The number $\eta(w_*)={\rm deg}(p(q))$ satisfies
 $\eta(w_*)=d_1+\cdots +d_l$. Moreover we have $\eta(w_*)=d$, the multiplicity
 of the singularity of the Painlev\'e divisor ${\mathcal D}_0=\cup_{j=1}^l{\mathcal D}_j$ at the point
 $p_o$,  for any semisimple Lie algebra not
 containing factors of type $E$ or $F$.
 We have the following formulas for  $\eta(w_*)$. The number $\eta(w_*)$ is, in each case,  the complex dimension of any Borel subalgebra of  ${\rm Lie}(K({\mathbb C}))$ $($and of  ${\rm Lie}(\check K({\mathbb C})))$,
 
 \begin{itemize}
\item[]  $A_l$:  $\eta(w_*) =\frac{l(l+2)}{4}$ if $l$ is even;$\quad\eta(w_*)=\frac{(l+1)^2}{4}$ if $l$ is odd,
\item[]  $B_l$  or $C_l$:  $\eta(w_*)=\frac{l(l+1)}{2}$,
\item[]  $D_l$:  $\eta(w_*)=\frac{l^2}{2}$ if $l$ is even;$\quad\eta(w_*)=\frac{l^2-1}{2}$ if $l$ is odd,
\item[]  $E_l$:  $\eta(w_*)=20$ if $l=6$;$\quad\eta(w_*)=35$ if $l=7$;$\quad\eta(w_*)=64$ if $l=8$,
\item[]  $F_4$:  $\eta(w_*)=14$,
\item[]  $G_2$:  $\eta(w_*)=4$.

\end{itemize}
 \end{Proposition}
 \begin{Proof} The explicit computation of $\eta(w_*)$  follows from  Theorem \ref{Kpq} which expresses $p(q)$ in terms of the order of certain Chevalley groups. Using  Proposition \ref{direct} for any semisimple Lie algebras with no $E$ or $F$ factors, we obtain  $\eta(w_*)=d$, the dimension of any Borel subalgebra of ${\rm Lie}(K({\mathbb C}))$. We also expect that $\eta(w_*)=d$ is true for all
 semisimple Lie algebra (Conjecture \ref{tangentcone}).
 \end{Proof}
  


\bibliographystyle{amsalpha}

\end{document}